\documentclass[11pt, reqno]{amsart}
\usepackage[T1]{fontenc}
\usepackage[utf8]{inputenc}
\usepackage{lmodern}
\usepackage{amsmath,amssymb,mathtools,mathrsfs}
\usepackage[margin=1.1in]{geometry}
\usepackage[expansion=false]{microtype}
\usepackage{enumitem,needspace,booktabs,tabularx}
\usepackage[all,cmtip]{xy}
\usepackage{xcolor}
\usepackage[hidelinks]{hyperref}
\hypersetup{colorlinks,	linkcolor={red!50!black},citecolor={blue!50!black},	urlcolor={blue!75!black},
	pdftitle={Arithmetic Kac theory for ramified species},pdfauthor={Jiarui Fei}}
\numberwithin{equation}{section}

\newtheorem{theorem}{Theorem}[section]
\newtheorem{proposition}[theorem]{Proposition}
\newtheorem{lemma}[theorem]{Lemma}
\newtheorem{corollary}[theorem]{Corollary}

\theoremstyle{definition}

\newtheorem{example}[theorem]{Example}

\theoremstyle{remark}
\newtheorem{remark}[theorem]{Remark}
\newtheorem*{theoremA}{Theorem A}
\newtheorem*{theoremB}{Theorem B}
\newtheorem*{theoremC}{Theorem C}

\newcommand{\op}[1]{\operatorname{#1}}
\newcommand{\mr}[1]{\mathrm{#1}}
\newcommand{\mb}[1]{\mathbb{#1}}
\newcommand{\mc}[1]{\mathcal{#1}}

\newcommand{\ms}[1]{\mathscr{#1}}
\newcommand{\OO}{\mc O}
\newcommand{\RR}{\mc R}
\newcommand{\GG}{\Gamma}
\newcommand{\Ind}{\op{Ind}}
\newcommand{\rk}{\underline{\op{rk}}}
\newcommand{\dimv}{\underline{\dim}}
\newcommand{\End}{\op{End}}
\newcommand{\Hom}{\op{Hom}}
\newcommand{\Ext}{\op{Ext}}
\newcommand{\rad}{\op{rad}}
\newcommand{\Gal}{\op{Gal}}
\newcommand{\Tr}{\op{Tr}}
\newcommand{\Lie}{\op{Lie}}
\newcommand{\Aut}{\op{Aut}}
\newcommand{\Br}{\op{Br}}
\newcommand{\lcm}{\op{lcm}}
\newcommand{\supp}{\op{supp}}
\newcommand{\Res}{\op{Res}}
\newcommand{\Rep}{\op{rep}}
\newcommand{\modfd}{\op{mod}_{\mr{fd}}}
\newcommand{\Qbar}{\widetilde Q}
\newcommand{\Kbar}{\overline K}
\newcommand{\iso}{\cong}
\newcommand{\into}{\hookrightarrow}
\newcommand{\unfoldcomp}{{c_{\Qbar}}}

\newcommand{\quiverdiagram}[1]{\vcenter{\hbox{$\displaystyle #1$}}}

\setlist[enumerate]{label=\textup{(\roman*)},leftmargin=2.3em,itemsep=2.9pt}

\title[Arithmetic Kac theory]{Arithmetic Kac theory for ramified species}
\author{Jiarui Fei}
\address{School of Mathematical Sciences, Shanghai Jiao Tong University, China}
\email{jiarui@sjtu.edu.cn}
\thanks{\hspace{-.19in} The author was supported in part by National Natural Science Foundation of China (Nos.~12131015 and 12571038)}

\date{}

\begin{document}
\begin{abstract}
Let $K=k((t))$, where $k$ is algebraically closed of characteristic
zero, and let $\Lambda$ be an acyclic ramified species over $K$.
An indecomposable $\Lambda$-module can split after extending $K$.
The number of its distinct geometric indecomposable summands
is the degree of its residue endomorphism field
$\End_\Lambda(M)/\rad\End_\Lambda(M)$.
We use the positive roots and cyclic automorphism of the unfolded
quiver to determine the number of isomorphism classes at every nonzero
rank and residue degree $e\ge1$: it is zero, an explicitly computable
finite number, or $|K|$.
For every invariant imaginary root of the unfolding, including
at nonprimitive ranks, we construct a smooth invariant family of
absolutely indecomposable representations whose rational quotient by
the change-of-bases group has dimension $1-q_d(\beta)$.
We characterize ranks admitting absolutely Schur representations
by the Schur roots of the unfolded quiver.

For the associated species over finite fields, we realize the polynomials
counting absolutely indecomposable representations as graded traces
of the unfolding automorphism on preprojective BPS cohomology.
Galois-orbit formulas and M\"obius inversion then recover the counts
at each residue degree. Finally, we prove an all-rank symmetrizable
analogue of Kac's constant-term conjecture: the constant term of the
polynomial counting all indecomposables, rather than only absolutely
indecomposable ones, equals the corresponding Kac--Moody root
multiplicity. This constant is independent of the symmetrizer and
modulation.
\end{abstract}
\maketitle

\begingroup\small
\tableofcontents
\endgroup
\clearpage

\section*{Introduction}

Kac's theorem identifies the dimension vectors of indecomposable quiver
representations over an algebraically closed field with positive roots.
Over other fields, indecomposables can split after scalar extension.
We study this phenomenon for species over Laurent-series fields,
determine the possible degrees of their residue endomorphism fields,
and relate finite-field counts to preprojective cohomology and
symmetrizable root multiplicities.

Gabriel proved that a connected acyclic quiver has finitely many
indecomposable isomorphism classes over an algebraically closed field
exactly when its underlying graph is Dynkin of type $A$, $D$, or $E$;
these classes correspond to positive roots \cite{Gabriel}. Kac extended
the correspondence to quivers without loops: real roots give unique
indecomposables and imaginary roots give infinitely many
\cite{Kac80,Kac83}.

Kac also introduced finite-field counts. A representation is
\emph{absolutely indecomposable} if it remains indecomposable over an
algebraic closure. For a quiver $Q$, the Kac polynomial
$A_{Q,\alpha}(T)$ counts such representations of dimension $\alpha$
over $\mb F_q$ at $T=q$ \cite{Kac83}. Kac conjectured coefficient
positivity and equality of the constant term with the corresponding
root multiplicity. The doubled-quiver representation space is the
cotangent bundle of the original representation space, and its
moment-map fibres describe deformed preprojective representations.
Crawley--Boevey characterized dimension vectors supporting simple
modules over these algebras \cite{CBmoment}. Crawley--Boevey and
Van den Bergh proved both conjectures for indivisible dimension vectors
using cohomology of smooth moduli spaces \cite{CBVdB}. Hausel proved
the constant-term conjecture in general \cite{Hausel}, and
Hausel--Letellier--Rodriguez-Villegas proved positivity \cite{HLRV}.

For symmetrizable Cartan matrices, \emph{species} assign division
algebras to vertices and bimodules to arrows. Dlab--Ringel developed
their finite and tame representation theory \cite{DlabRingel}.
Alternatively, one may start with an ordinary quiver and an
\emph{admissible automorphism}, whose vertex orbits contain no arrows.
Hubery established the root correspondence for sums of indecomposables
along such orbits, including uniqueness at real roots \cite{Hubery}.

Over a finite field, a quiver automorphism can be combined with
Frobenius. Deng and Du established the resulting correspondence between
representations of the fixed algebra and representations with compatible
Frobenius structures \cite{DengDu}. Their Theorem~5.1 and
Corollary~5.2 identify indecomposables with Frobenius orbits and identify
the orbit length with the degree of the endomorphism algebra modulo its
radical. They also describe scalar extension and prove polynomiality and
orientation independence of the total indecomposable count
\cite[Corollaries~6.4 and~9.2]{DengDu}. Their counting argument uses an
intertwiner formula of Hua and Lin \cite{HuaLin}. Deng and Xiao prove
the finite-field species root correspondence by Ringel--Hall methods
\cite{DengXiao}; see also \cite[Section~6]{Hubery}.

Gei\ss--Leclerc--Schr\"oer realized symmetrizable Cartan data by
algebras with truncated-polynomial vertex rings \cite{GLSI,GLSII}.
Modules free over each vertex ring are \emph{locally free}, and their
vertex ranks replace dimension vectors. Their paper \cite{GLSrigid} 
compares the GLS algebra $H(C,D,\Omega)$, an integral order over a
power-series ring, and the corresponding generic hereditary species
through reduction and localization. This comparison concerns rigid modules,
which have no self-extensions. The Laurent-series species studied here
arise from that construction. Galois descent for regular representations
of such species also appears in \cite{GR,GRR}.

The cohomological approach is based on the work of Davison and
Meinhardt. Their cohomological Poincar\'e--Birkhoff--Witt (PBW) theorem
expresses representation-stack cohomology as a symmetric algebra on
cohomology of \emph{BPS sheaves}, together with a scalar factor
\cite{DM}. Davison proved purity for preprojective stacks in every
dimension and identified the associated BPS cohomology with Kac
polynomials \cite{DavisonStacks}. He also identified the degree-zero
part of the preprojective BPS Lie algebra with the negative half of the
ordinary Kac--Moody algebra \cite[Theorem~6.6]{DavisonBPS}.
We apply these constructions with the automorphism of the unfolding.
The resulting traces count
representations of the species. The constant-term comparison also uses
the Verma-module twining theorem of Fuchs--Schellekens--Schweigert
\cite{FSS}.

\subsection*{The species and its representations}
Let $k$ be algebraically closed of characteristic zero, and choose
compatible roots of $t$ in an algebraic closure $\Kbar$ of
\[ K=k((t)),\qquad K_m=k((t^{1/m}))\subset\Kbar\quad(m\ge1). \]
Let $I$ be a finite vertex set, let $d=(d_i)_{i\in I}$ be a tuple of
positive integers, and put
\[ R=\prod_{i\in I}K_{d_i},\qquad D=\op{diag}(d_i). \]
Write $e_i$ for the idempotent that is one in the $i$th factor and zero
elsewhere. Let $A$ be a finite-dimensional $R$-bimodule whose left and
right $K$-actions agree, and set $A_{ij}=e_iAe_j$. We assume that
$A_{ii}=0$ and that the directed graph with $i\to j$ whenever
$A_{ij}\ne0$ is acyclic. The vertex fields and arrow bimodules form an
\emph{acyclic ramified species}; their choice, including the two field
actions on each bimodule, is its \emph{modulation}. Its algebra is
\[ \Lambda=T_R(A)=R\oplus A\oplus(A\otimes_R A)\oplus\cdots . \]
Acyclicity makes this sum finite. We use right modules and multiply
paths from left to right.

Throughout the Laurent-series part of the paper, a representation $M$
is a right $\Lambda$-module finite-dimensional \emph{over $K$}.
Equivalently, it consists of finite-dimensional $K_{d_i}$-spaces
$M_i=Me_i$ and $K_{d_j}$-linear maps
\[ M_i\otimes_{K_{d_i}}A_{ij}\longrightarrow M_j. \]
Its \emph{rank vector} measures the dimensions over the vertex fields:
\[ \rk M=(\dim_{K_{d_i}}M_i)_{i\in I}=\beta,\qquad \dim_K M=\sum_i d_i\beta_i. \]
For an extension $E/K$, write $\Lambda_E=E\otimes_K\Lambda$ and
$M_E=E\otimes_K M$. A \emph{geometric representation} is a module over
$\Lambda_{\Kbar}$.

The Cartan matrix of the species is $C=(c_{ij})$, where
\[ c_{ii}=2,\qquad c_{ij}=-\frac{\dim_K A_{ij}+\dim_K A_{ji}}{d_i}\quad(i\ne j). \]
Then $DC$ is symmetric, so $D$ is a \emph{symmetrizer} of $C$.
Every symmetrizable Cartan matrix with an integral symmetrizer and an
acyclic orientation occurs in this construction
(Proposition~\ref{prop:standard}). The pair $(C,D)$ need not determine
the modulation.

We use $\mb N=\{0,1,2,\ldots\}$ and
$\supp\beta=\{i\in I:\beta_i\ne0\}$. For $\beta\in\mb N^I$, put
\[ q_d(\beta)=\tfrac12\beta^{\mathsf T}DC\beta,\qquad p_d(\beta)=1-q_d(\beta),\qquad g_d(\beta)=\gcd(d_i\beta_i:i\in I), \]
with $g_d(0)=0$. A nonzero vector is \emph{primitive} if its coordinates
have greatest common divisor one. We write $\Delta^+(C)$ for the
positive roots and $\mathfrak g(C)_\beta$ for the corresponding root
space of the complex Kac--Moody algebra, and set
$\mathfrak g(C)_\beta=0$ when $\beta$ is not a root. A root is \emph{real} if it is a Weyl translate of a
simple root, and \emph{imaginary} otherwise. An imaginary root of
quadratic value zero is \emph{isotropic}.

\subsection*{Galois orbits and residue endomorphism fields}
For an indecomposable $M$, the quotient
$\End_\Lambda(M)/\rad\End_\Lambda(M)$ is a division algebra. Over $K$
it is a finite extension $K_e$ (Lemma~\ref{lem:fields}). We call it the
\emph{residue endomorphism field} and call $e$ its \emph{residue degree}.
This degree concerns the endomorphism algebra, not the valuation residue field of $K_e$, which
is always $k$.

The degree $e$ has a representation-theoretic interpretation:
$M_{\Kbar}$ is the direct sum of exactly $e$ pairwise nonisomorphic
indecomposables, each occurring once
(Proposition~\ref{prop:geometric-orbit}). In particular, $e=1$ is
equivalent to absolute indecomposability. For every $\beta\in\mb N^I$
and $e\ge1$, define
\[ \Ind_{\beta,e}(\Lambda)=\{[M]:M\text{ indecomposable},\ \rk M=\beta,\ \End_\Lambda(M)/\rad\End_\Lambda(M)\iso K_e\}. \]
Nonemptiness at a nonzero rank will imply $e\mid g_d(\beta)$.

For example, the vertex simple $S_i$ has $K_{d_i}$ at vertex $i$ and
zero elsewhere. Its rank is $\epsilon_i$, where $\epsilon_i$ is the
$i$th standard basis vector, and its endomorphism field is $K_{d_i}$.
After scalar extension it becomes a sum of $d_i$ vertex simples.
Thus even a simple representation of a species need not be absolutely
indecomposable.

To describe the general decomposition, put $N=\lcm_i d_i$. A splitting
\[ \Lambda_{K_N}\iso K_N\Qbar \]
identifies the species with an ordinary acyclic quiver after extending
the field (Proposition~\ref{prop:unfold}). Choose a generator $\sigma$
of $\Gal(K_N/K)$. In suitable arrow bases its action is
\[ \tau\!\left(\sum_p c_p p\right)=\sum_p\sigma(c_p)\,a(p), \]
where $p$ ranges over paths and $a$ is a quiver automorphism. Thus the
descent acts on coefficients \emph{and} on the quiver. It is
$\sigma$-semilinear, in the precise sense that
$\tau(cx)=\sigma(c)\tau(x)$. We call $\Qbar$ the \emph{unfolding}.

The vertices above $i$ are $(i,r)$, with
$r\in\mb Z/d_i\mb Z$, and $a(i,r)=(i,r+1)$. No arrow joins
two vertices above the same $i$, so $a$ is admissible. Write $\Qbar_0$ for its vertex set,
$\widetilde C$ for its symmetric Cartan matrix, and
$\Delta^+(\Qbar)=\Delta^+(\widetilde C)$. On dimension vectors,
$(a\alpha)_v=\alpha_{a^{-1}v}$. The \emph{repetition map} is
\[ \iota:\mb Z^I\longrightarrow\mb Z^{\Qbar_0},\qquad (\iota\beta)_{(i,r)}=\beta_i. \]
It identifies rank vectors with $a$-invariant dimension vectors:
$\dimv M_{\Kbar}=\iota\beta$ when $\rk M=\beta$.

For $e\ge1$, define the \emph{orbit norm}
\[ \mc N_e^a\alpha=\sum_{j=0}^{e-1}a^j\alpha. \]
If $a^e\alpha=\alpha$, this sum is $a$-invariant. The geometric
dimensions of a representation in $\Ind_{\beta,e}(\Lambda)$ therefore
give orbits in the finite set
\[ \ms R_{\beta,e}=\{\alpha\in\Delta^+(\Qbar):a^e\alpha=\alpha,\ \mc N_e^a\alpha=\iota\beta\}/\langle a\rangle. \]
This set is finite, since $0\le\alpha\le\iota\beta$ coordinatewise.
The \emph{dimension period} $s(\alpha)$ is the least positive integer
with $a^{s(\alpha)}\alpha=\alpha$. It may be smaller than the
\emph{isomorphism period}, the size of a geometric representation's
Galois orbit. Write $\ms R^{\mr{im}}_{\beta,e}$ for the
imaginary root orbits, and $\ms R^{\mr{re,exact}}_{\beta,e}$
for the real root orbits with $s(\alpha)=e$.

For real roots, the dimension period must be exactly $e$;
Example~\ref{ex:doubled-indefinite} shows why this condition matters.

\begin{theoremA}[Arithmetic classification]
For every nonzero rank $\beta$ and every $e\ge1$,
\[
\#\Ind_{\beta,e}(\Lambda)=
\begin{cases}
|K|,&\ms R^{\mr{im}}_{\beta,e}\ne\varnothing,\\[2pt]
|\ms R^{\mr{re,exact}}_{\beta,e}|,
 &\ms R^{\mr{im}}_{\beta,e}=\varnothing.
\end{cases}
\]
The criterion is independent of the acyclic orientation of the chosen
modulation. The set is empty unless $e\mid g_d(\beta)$.
\end{theoremA}

\noindent\emph{The geometry behind Theorem A.}
For $e=1$, Theorem~A says that absolutely indecomposable representations
exist precisely when $\iota\beta$ is a positive root; real roots give a
unique class. Its proof uses a stronger geometric assertion for
imaginary roots.

Let $V_\beta$ be the affine space of representations
on fixed vertex spaces $K_{d_i}^{\beta_i}$, and let $G_\beta$ be the
algebraic group of changes of these bases. Explicitly,
\[ G_\beta=\prod_i\Res_{K_{d_i}/K}\mr{GL}_{\beta_i},\qquad G_\beta(K)=\prod_i\mr{GL}_{\beta_i}(K_{d_i}). \]
Here $\Res$ denotes Weil restriction.
Set $\lVert\beta\rVert_d^2:=\sum_i d_i\beta_i^2=\dim G_\beta$. For every imaginary root $\iota\beta$, Theorem~\ref{thm:absolute-root}
gives a smooth geometrically integral, $G_\beta$-invariant locally closed
subvariety $W\subseteq V_\beta$, with a $K$-point, such that every
geometric point is indecomposable. The endomorphism dimension is constant
on $W$; writing $h$ for this dimension, we have
\[ \nu(W)=\dim W-\lVert\beta\rVert_d^2+h=1-q_d(\beta)>0. \]
Each geometric orbit has dimension $\dim G_\beta-h$, so a rational
quotient of $W$ has dimension $1-q_d(\beta)$. This is also the degree
of the ordinary Kac polynomial $A_{\Qbar,\iota\beta}(T)$.
No primitivity assumption is needed. The same assertions hold after
every finite extension of $K$, for the corresponding partially split
species.

When $q_d(\beta)<0$, deformed preprojective representations let us
change orientation without losing parameters. For isotropic roots we
descend a rational curve of isomorphism classes to $K$. Varying a rational invariant
over extensions prescribes Galois periods beyond the dimension period, leading from Theorem~\ref{thm:absolute-root} to
Theorem~\ref{thm:allrank-strata}. At real roots, geometric uniqueness
forces the two periods to agree, so shorter-period roots contribute no
class.

Theorem~\ref{thm:allrank-strata} proves Theorem~A. It gives a finite
root calculation in every Cartan type, rather than normal forms for
individual representations in a wild family. As a consequence, the
indecomposable ranks are exactly $\Delta^+(C)$; real roots give unique
classes and imaginary roots give $|K|$ classes
(Corollary~\ref{cor:full-kac}). In affine type,
Proposition~\ref{prop:affine-strata} also describes the individual
classes through homogeneous and exceptional tubes. Examples
\ref{ex:two-modulations} and~\ref{ex:14-affine} show why their residue
degrees carry information not determined by the total number of classes.

\subsection*{Generic Schur representations}
An indecomposable representation is \emph{Schur} if its endomorphism
algebra is a division algebra; it is \emph{absolutely Schur} if that
algebra is $K$. Proposition~\ref{prop:generic-schur-summands} gives,
at every rank $\beta$, a nonempty open subset $U_\beta\subseteq V_\beta$
defined over $K$ such that every indecomposable summand of every
$M\in U_\beta(K)$ has endomorphism algebra $K_e$ for some $e$.
Corollary~\ref{cor:absolute-schur-ranks} characterizes the nonempty
absolutely Schur loci: they occur exactly when $\iota\beta$ is a
Schur root of the unfolding, or equivalently when a general geometric
representation of rank $\beta$ is indecomposable. Existence of an
arithmetic Schur representation alone does not imply this geometric
genericity; see Remark~\ref{rem:arithmetic-schur-genericity}.

\subsection*{Why we count all indecomposables}
The same unfolding defines a species over a finite field:
\[ \Lambda_{q,a}=(\overline{\mb F}_q\Qbar)^{a\circ\mr{Frob}_q},\qquad \mr{Frob}_q(x)=x^q. \]
Its vertex fields are $\mb F_{q^{d_i}}$. In this setting a
representation is a finite-dimensional right $\Lambda_{q,a}$-module
over $\mb F_q$, and rank is measured over these vertex fields.

Let $I_{\beta,e}(q)$ count its indecomposable isomorphism classes of
rank $\beta$ with residue endomorphism field $\mb F_{q^e}$. For
$\beta\ne0$, set
\[ I_{\beta,a}(q)=\sum_{e\mid g_d(\beta)}I_{\beta,e}(q). \]
Thus $I_{\beta,1}$ counts absolute indecomposables and $I_{\beta,a}$
counts all indecomposables. The dependence on $a$ is suppressed in
$I_{\beta,e}$.

For the split quiver, $I_{\beta,1}=A_{\Qbar,\beta}$, but the total
count can differ from the absolute count even when $a=\mr{id}$.
For nontrivial $a$, the Frobenius descent also changes the absolute count.
The vertex simple gives an even more immediate reason to use the total
count in a folded root-multiplicity formula:
\[ I_{\epsilon_i,d_i}(T)=1,\qquad I_{\epsilon_i,a}(T)=1,\qquad I_{\epsilon_i,1}(T)=0\quad(d_i>1). \]
Its rank is a simple root of multiplicity one. An absolute count would
omit precisely this contribution. We therefore use $A$ for ordinary
Kac polynomials and $I$ for the counts refined by residue degree and
their sum.

\subsection*{Cohomological traces and root multiplicities}
We now work with the ordinary unfolding over $\mb C$ and use
rational cohomology. The BPS sheaves for its preprojective algebra give
finite-dimensional rational multiplicity spaces $B_{\alpha,j}$,
with only finitely many nonzero $j$ for each fixed $\alpha$. In the dual convention specified in
Subsection~\ref{subsec:bps-objects},
\[ \sum_{j\ge0}\dim_{\mb Q}B_{\alpha,j}\,T^j=A_{\Qbar,\alpha}(T). \]
This is an identity of dimensions. The sheaf construction also supplies
maps $B_{\alpha,j}\to B_{u\alpha,j}$ for quiver automorphisms $u$.
For a power $u$ of $a$ fixing $\alpha$, define the \emph{trace polynomial}
\[ P_{\alpha,u}(T)=\sum_{j\ge0}\Tr(u\mid B_{\alpha,j})T^j; \]
set it to zero when $u\alpha\ne\alpha$. In particular,
$P_{\alpha,\mr{id}}=A_{\Qbar,\alpha}$.

\begin{theoremB}[Arithmetic trace formula]
For every nonzero rank $\beta$ and every finite field $\mb F_q$,
\[ I_{\beta,1}(q)=P_{\iota\beta,a}(q). \]
The formula holds for every power of $a$, using the vertex orbits and
rank coordinates of that power. Every residue-degree count is the
evaluation of the polynomial
\[ I_{\beta,e}(T)=\frac1e\sum_{\substack{\alpha\in\mb N^{\Qbar_0}\setminus\{0\}:\ a^e\alpha=\alpha\\\mc N_e^a\alpha=\iota\beta}}\ \sum_{r\mid e}\mu_{\mr{ar}}(e/r)P_{\alpha,a^r}(T^r). \]
Here $\mu_{\mr{ar}}$ is the arithmetic M\"obius function. The
absolute polynomial belongs to $\mb Z[T]$. The other polynomials
belong to $\mb Q[T]$ and take integer values at prime powers, but
need not have integral coefficients.
\end{theoremB}

Theorem~B refines Deng--Du's polynomiality of the total count
\cite[Corollary~9.2]{DengDu}. The divisor sum counts Galois
orbits of length $e$. After extension to $\mb F_{q^r}$, descent uses
$a^r\circ\mr{Frob}_{q^r}$, so vertex orbits can split.
The proof compares Jordan-form counts with equivariant PBW traces.
Crucially, symmetric powers are formed in the full dimension-vector
grading before restricting to invariant degrees: an orbit of moving
dimensions may have invariant total dimension. The same comparison
gives the recursion \eqref{eq:primitive-recursion}.

When $g_d(\beta)=1$, the trace polynomial also has a realization on a
smooth stable preprojective moduli space
(Corollary~\ref{cor:smooth-trace}).

Section~\ref{sec:polynomial-strata} relates these polynomials to
the three cardinality alternatives in Theorem~A and determines
their degrees and leading coefficients. At primitive rank,
$I_{\beta,e}(T)\in\mb Z[T^e]$ (Corollary~\ref{cor:primitive-specialization}).

The following result is an all-rank symmetrizable analogue of
Kac's constant-term conjecture. The distinction in what is counted
is essential: the classical conjecture concerns absolutely indecomposable
representations, whereas the folded statement uses the total number of
indecomposable species representations.

\begin{theoremC}[Folded constant term]
For every nonzero rank $\beta$,
\[ I_{\beta,a}(0)=\dim_{\mb C}\mathfrak g(C)_\beta, \]
with the right side zero at a nonroot. This total constant is
independent of the integral symmetrizer and the modulation.
\end{theoremC}

The proof identifies the degree-zero BPS action with the action on
Chevalley generators and then applies the Verma-module twining theorem
(Theorem~\ref{thm:folded-constant}). Unlike ordinary Kac polynomials,
the arithmetic counting polynomials need not have nonnegative
coefficients; thus the theorem is a constant-term analogue, not
a coefficient-positivity statement.

Signed coefficients reflect automorphism traces, rather than
negative dimensions of cohomology spaces
(Examples~\ref{ex:arithmeticfamily} and~\ref{ex:indefinite-four-degrees}).

Finally, Section~\ref{sec:orders} compares integral GLS orders with
their generic-fibre ramified species. We show that localization
need not preserve indecomposability and that lattice lifts of the
nonroot-rank locally free modules of Lin--Su \cite{LinSu} have decomposable
generic fibres. We also give a connected, minimal-symmetrizer example
contradicting two GLS conjectures concerning irreducibility of
representation varieties and generic endomorphism algebras at Schur roots.

\subsection*{Conventions and organization}
For an ordinary quiver $Q$, its vertex and arrow sets are $Q_0$ and
$Q_1$. We write $\Delta^+(Q)$ for the roots of its symmetric Cartan
matrix $C_Q$, and put
\[ q_Q(\alpha)=\tfrac12\alpha^{\mathsf T}C_Q\alpha,\qquad p_Q(\alpha)=1-q_Q(\alpha),\qquad |\alpha|=\sum_{v\in Q_0}\alpha_v. \]
The symmetric form is $(\alpha,\eta)_Q=\alpha^{\mathsf T}C_Q\eta$;
we omit $Q$ when it is fixed. In particular,
$q_{\Qbar}(\iota\beta)=q_d(\beta)$. For the folded root system,
$s_i\beta=\beta-(C\beta)_i\epsilon_i$ is the simple reflection.
The notation $0<\gamma<\alpha$ means coordinatewise containment with
$\gamma\ne0,\alpha$, not strict inequality at every coordinate.

Over a field $F$, $V_\alpha=\Rep(Q,\alpha)$ is the affine space of
arrow maps on the spaces $F^{\alpha_v}$, and
$G_\alpha=\prod_v\mr{GL}_{\alpha_v}$ is its change-of-bases group.
For species these constructions are given in Section~\ref{sec:moment}.
We also use $\Res$ for restriction of scalars on modules; its
argument distinguishes this from Weil restriction on groups.

Sections~\ref{sec:unfold}--\ref{sec:descent} develop unfolding,
Galois descent, and the generic Schur results.
Sections~\ref{sec:moment}--\ref{sec:classification} construct the
families of absolutely indecomposable representations and use them to
prove Theorem~A; the affine applications follow in
Section~\ref{sec:affine-strata}.
Sections~\ref{sec:count}--\ref{sec:polynomial-strata} prove
Theorems~B and~C and compare the counting polynomials with the sets
classified in Theorem~A. The proof of Theorem~A does not use BPS
cohomology. Section~\ref{sec:orders} treats integral orders,
localization, and counterexamples to two GLS conjectures.

\clearpage

\section{Ordinary species and cyclic unfoldings}\label{sec:unfold}

In species diagrams, arrows are labelled by bimodules.
Write
\[ h_{ij}=\gcd(d_i,d_j),\qquad \ell_{ij}=\lcm(d_i,d_j), \qquad N=\lcm_{i\in I}d_i. \]
The fields $K_{d_i}$ and $K_{d_j}$ have intersection $K_{h_{ij}}$ and
compositum $K_{\ell_{ij}}$. Since the extensions are separable,
\begin{equation}\label{eq:tensorfields} K_{d_i}\otimes_K K_{d_j}\iso \prod_{s\in\mb Z/h_{ij}\mb Z}K_{\ell_{ij}}.\end{equation}
In particular, each $A_{ij}$ is a direct sum of simple bimodules of
$K$-dimension $\ell_{ij}$. Its dimension is divisible by both $d_i$
and $d_j$, as used in the definition of $C$.

\begin{proposition}\label{prop:standard}
Let $C$ be a symmetrizable generalized Cartan matrix, let
$D=\op{diag}(d_i)$ be an integral symmetrizer, and choose an
acyclic orientation $\Omega$. For $i\to j$ in $\Omega$, set
\[ p_{ij}=\gcd(-c_{ij},-c_{ji}),\qquad A_{ij}=K_{\ell_{ij}}^{\,p_{ij}}, \]
with the natural actions of its two subfields. Set all other arrow
bimodules equal to zero. The Cartan matrix of this species is $C$.
\end{proposition}
\begin{proof}
The equality $d_i(-c_{ij})=d_j(-c_{ji})$ implies
$-c_{ij}=p_{ij}\frac{d_j}{h_{ij}}$ and $-c_{ji}=p_{ij}\frac{d_i}{h_{ij}}$.
Hence $\dim_K A_{ij}=p_{ij}\ell_{ij}=d_i(-c_{ij})$, which proves
the assertion.
\end{proof}

\begin{proposition}\label{prop:unfold}
There is a finite acyclic quiver $\Qbar$, an admissible automorphism
$a$ of $\Qbar$ of order dividing $N$, and an isomorphism
$K_N\otimes_K\Lambda\iso K_N\Qbar$
with the following properties.
\begin{enumerate}
\item The vertices over $i$ are $(i,r)$, with
$r\in\mb Z/d_i\mb Z$, and $a(i,r)=(i,r+1)$.
\item Under the isomorphism, a generator of $\Gal(K_N/K)$ acts by
coefficient conjugation together with $a$ on vertices and arrows.
\item The folded Cartan matrix is $C$, and the vertex-orbit sizes
are the diagonal entries of $D$.
\item For every $\Lambda$-module $M$,
\begin{equation}\label{eq:rank-dim} \dimv(M_{\Kbar})=\iota(\rk M),\qquad (\iota\beta)_{(i,r)}=\beta_i.\end{equation}
\end{enumerate}
Here admissibility means that no arrow joins two vertices in the same
$a$-orbit.
\end{proposition}
\begin{proof}
We split each vertex field as
$K_N\otimes_K K_{d_i}\iso\prod_{r\in\mb Z/d_i\mb Z}K_N$,
and the Galois generator permutes the factors cyclically. By
\eqref{eq:tensorfields}, it is enough to split one simple arrow
bimodule. Such a bimodule contributes one arrow between each pair of
sheets $(i,r),(j,s)$ satisfying
\[ r-s\equiv c\pmod {h_{ij}} \]
for a suitable fixed residue $c$. There are $\ell_{ij}$ such pairs. Simultaneous translation $(r,s)\mapsto(r+1,s+1)$ permutes them
transitively. A direct sum decomposition of the bimodule gives the
same construction with multiplicities. We choose the arrow basis from the primitive idempotents of these
split field factors. In this basis the Galois generator permutes the
arrows without additional scalar factors.

Tensor algebras commute with scalar extension. Thus the split vertex
algebra and the split arrow bimodule give $K_N\Qbar$. A directed cycle
in $\Qbar$ would project to a directed cycle in the support of $A$,
so $\Qbar$ is acyclic. The absence of $A_{ii}$ proves admissibility.

There are $d_i$ vertices in the orbit over $i$, and the total number
of edges between the orbits over $i$ and $j$ is
$\dim_K A_{ij}+\dim_K A_{ji}$. The folded Cartan entry is therefore
$c_{ij}$. Finally, a $K_{d_i}$-space of dimension $r_i$ splits into
$d_i$ spaces of dimension $r_i$, proving \eqref{eq:rank-dim}.
\end{proof}

For the standard species, this agrees with the unfolding in
\cite[Lemma~3.3]{Pfeifer}. The preceding proof permits all the factors
in \eqref{eq:tensorfields}, including the factors in which the two
embeddings of $K_{h_{ij}}$ differ. The corresponding finite-field
construction is \cite[Proposition~6.3]{DengDu}.

For a quiver automorphism $u$, write ${}^uY$ for the representation
obtained by relabelling vertices and arrows by $u$; thus
$\dimv({}^uY)=u\dimv Y$. We shall use the following consequences of
Hubery's folding theorem. Let $\widetilde C$ be
the symmetric Cartan matrix of $\Qbar$, and let
\[ f:(\mb Z^{\Qbar_0})^{\langle a\rangle}\longrightarrow\mb Z^I \]
be the inverse of $\iota$.
For a root $\alpha$ of $\widetilde C$, let $s$ be its least positive
period under $a$. Its orbit norm is $\mc N_s^a\alpha$.
Hubery proves the following facts
\cite[Proposition~4 and Theorem~1]{Hubery}:
\begin{enumerate}
\item $f(\mc N_s^a\alpha)$ is a root of $C$. If it is real,
$\alpha$ is real and is determined up to $a$-orbit.
\item Every positive root $\beta$ of $C$ is $f(\dimv X)$ for an
orbit sum
\[ X=Y\oplus{}^aY\oplus\cdots\oplus{}^{a^{s-1}}Y \]
of pairwise nonisomorphic indecomposable $k\Qbar$-modules, where $s$
is the least positive period of the isomorphism class of $Y$.
\item If $\beta$ is real, this orbit sum is unique up to isomorphism
and has $q_d(\beta)$ summands.
\end{enumerate}
We shall also use the \emph{fundamental region}
\[ \mc F(C)=\{\gamma\in\mb N^I\setminus\{0\}:\supp\gamma\text{ is connected},\ (C\gamma)_i\leq0\text{ for every }i\}. \]
The positive imaginary roots are the Weyl translates of vectors in
$\mc F(C)$ \cite[Chapter~5]{KacBook}. This description shows, in
particular, that every positive integral multiple of a positive
imaginary root is again a positive imaginary root.

\section{Galois descent and possible ranks}\label{sec:descent}

The descent arguments below apply to every finite-dimensional
$K$-algebra $B$. All $B$-modules in this section are finite-dimensional
over $K$. For a field extension $E/K$, we use
$B_E=E\otimes_K B$ and $M_E=E\otimes_K M$. We take $B=\Lambda$ in Subsection~\ref{sec:root}. Write
$\GG=\Gal(\Kbar/K)$. Its action on $B_{\Kbar}$ fixes $B$ and
conjugates coefficients, and twisting a module by this action gives
the Galois action on geometric isomorphism classes. We use $\Br(F)$
for the Brauer group of a field $F$; its vanishing means that every
central simple $F$-algebra is a matrix algebra over $F$.

The finite-field version of the orbit correspondence below, including
the residue endomorphism field, is due to Deng--Du
\cite[Theorem~5.1 and Corollary~5.2]{DengDu}. Over $K$ we use Brauer
vanishing over every finite extension in place of the commutativity
of finite division rings. The remaining arguments use standard restriction-of-scalars methods.

\begin{lemma}\label{lem:fields}
Every finite extension of $K$ inside $\Kbar$ is $K_m$ for a unique
positive integer $m$, and $\Gal(K_m/K)$ is cyclic of order $m$.
For every such extension $F$, one has $\Br(F)=0$. Consequently,
every finite-dimensional division algebra over $K$ is a commutative
finite extension of $K$.
\end{lemma}
\begin{proof}
The Newton--Puiseux theorem gives $\Kbar=\bigcup_{m\ge1}K_m$.
Since $k$ contains every root of unity, $K_m/K$ is cyclic of degree
$m$, and its subfields give the asserted list of finite extensions.

To show that the Brauer group vanishes, we prove that norms are
surjective from every finite extension. Write $F=k((s))$ and $E=F(s^{1/m})$. If
$u\in k[[s]]^\times$, Hensel's lemma gives $v\in F$ with $v^m=u$,
so $N_{E/F}(v)=u$. Moreover
$N_{E/F}(s^{1/m})=(-1)^{m-1}s$. Thus the norm contains every unit
and a uniformizer. A central simple $F$-algebra splits over a finite
cyclic extension, whose relative Brauer group is
$F^\times/N_{E/F}(E^\times)=0$; see \cite[Chapters~4 and~6]{GS}.
Finally the centre of a finite-dimensional division algebra over $K$
is such a field $F$, and $\Br(F)=0$ forces the division algebra to
be $F$ itself.
\end{proof}

\begin{lemma}\label{lem:basechange}
Let $B$ be a finite-dimensional $K$-algebra and let $M,N$ be
finite-dimensional $B$-modules. For every field extension $E/K$, we have
\[ E\otimes_K\Hom_B(M,N)\iso\Hom_{B_E}(M_E,N_E). \]
If $T$ is a finite-dimensional $K$-algebra, then
\[ \rad(T\otimes_K\Kbar)=(\rad T)\otimes_K\Kbar. \]
\end{lemma}
\begin{proof}
We express $B$-linearity as finitely many homogeneous linear equations
in the entries of a $K$-linear map. Extension of scalars preserves
their kernel. For the second assertion,
$(\rad T)_{\Kbar}$ is nilpotent. The semisimple algebra $T/\rad T$
is separable over the characteristic-zero field $K$, so its scalar
extension is semisimple. These two observations give the equality.
\end{proof}

\begin{lemma}\label{lem:ND}
If $M_{\Kbar}\iso N_{\Kbar}$, then $M\iso N$.
\end{lemma}
\begin{proof}
An isomorphism and its inverse are given by finitely many matrix entries,
so they are defined over a finite extension $E/K$. Restricting scalars
from $E$ gives $M^{[E:K]}\iso N^{[E:K]}$. Krull--Schmidt cancellation
implies $M\iso N$; compare \cite[Lemma~8.4.1]{DWbook}.
\end{proof}

\begin{proposition}\label{prop:geometric-orbit}
Let $M$ be an indecomposable finite-dimensional $B$-module. There is
a positive integer $e$ such that
\[ \End_B(M)/\rad\End_B(M)\iso K_e. \]
The module $M_{\Kbar}$ is a direct sum of $e$ pairwise nonisomorphic
indecomposables, each occurring with multiplicity one, and their
isomorphism classes form a single $\GG$-orbit.
\end{proposition}
\begin{proof}
The endomorphism algebra $\End_B(M)$ is local. Its quotient by the radical
is a division algebra, hence is isomorphic to $K_e$ for some $e$ by
Lemma~\ref{lem:fields}. Lemma~\ref{lem:basechange} gives
\begin{equation}\label{eq:semisimpleend} \End_{B_{\Kbar}}(M_{\Kbar})/\rad\End_{B_{\Kbar}}(M_{\Kbar}) \iso\Kbar\otimes_KK_e\iso\Kbar^{\,e}.\end{equation}
On the other hand, if the Krull--Schmidt decomposition is
$M_{\Kbar}=\bigoplus_j X_j^{m_j}$, its semisimple endomorphism algebra
is $\prod_j\op{Mat}_{m_j}(\Kbar)$. Comparison with
\eqref{eq:semisimpleend} gives $m_j=1$ and exactly $e$ summands.
The action of $\GG$ on the factors of $\Kbar\otimes_KK_e$ is
transitive, since these factors correspond to the embeddings of the
field $K_e$. These factors correspond to the isomorphism classes of the $X_j$.
The action on those classes is therefore transitive.
\end{proof}

\begin{proposition}[Orbit descent]\label{prop:descent}
Scalar extension induces a bijection
\[ \Ind(B)\ \longleftrightarrow\ \Ind(B_{\Kbar})/\GG, \]
where $\Ind$ denotes isomorphism classes of nonzero finite-dimensional
indecomposable modules. Under this bijection, an orbit corresponds to
the module whose geometric scalar extension contains each member of
the orbit exactly once.
\end{proposition}
\begin{proof}
Proposition~\ref{prop:geometric-orbit} gives the map, and
Lemma~\ref{lem:ND} proves injectivity.

Let $X$ be an indecomposable $B_{\Kbar}$-module. Its finitely many
matrix coefficients lie in a finite extension $E/K$, so it has a
model $X_E$ over $B_E$. This model is absolutely indecomposable.
Form the finite-dimensional $B$-module $V=\Res_{E/K}X_E$.
Since $E/K$ is separable,
\[ V_{\Kbar}\iso \bigoplus_{\sigma:E\into\Kbar}(X_E\otimes_{E,\sigma}\Kbar). \]
All the summands on the right belong to the single Galois orbit of
$X$, possibly with repetitions. Take any indecomposable summand $M$
of $V$. By Proposition~\ref{prop:geometric-orbit}, the summands of
$M_{\Kbar}$ form a nonempty Galois-stable subset of that orbit, with no repetitions.
Since the orbit is transitive, this subset is the whole orbit. Hence
$M$ has the required geometric decomposition.
\end{proof}

\begin{proposition}\label{prop:finitebase}
Let $M$ be indecomposable and let
$\End(M)/\rad\End(M)\iso K_e$. For $m\geq1$, put
$h=\gcd(e,m)$ and $\ell=\lcm(e,m)$. Then $M_{K_m}$ is a direct sum
of $h$ pairwise nonisomorphic indecomposables, each occurring once,
and each with residue endomorphism field $K_\ell$ over $K_m$.
In particular,
\[ M_{K_m}\text{ is indecomposable}\quad\Longleftrightarrow\quad \gcd(e,m)=1. \]
All geometric summands are defined over $K_m$ exactly when $e\mid m$.
\end{proposition}
\begin{proof}
We compute the semisimple endomorphism algebra after base change:
\[ \End(M_{K_m})/\rad\End(M_{K_m}) \iso K_m\otimes_KK_e \iso\prod_{\nu=1}^{h}K_\ell. \]
Primitive idempotents lift through the nilpotent radical. Each corner
is local, and there are no matrix factors, so the resulting $h$
indecomposables are pairwise nonisomorphic and occur once. Their
residue degree over $K_m$ is $\ell/m=e/h$. They are all absolutely
indecomposable exactly when $e/h=1$, or equivalently $e\mid m$.
\end{proof}

\begin{proposition}\label{prop:divisibility}
If $M$ is an indecomposable $\Lambda$-module of rank $\beta$ with
residue endomorphism field $K_e$, then
\[ e\mid d_i\beta_i\quad(i\in I),\qquad e\mid g_d(\beta). \]
Consequently, if $g_d(\beta)=1$, every indecomposable of rank $\beta$
is absolutely indecomposable.
\end{proposition}
\begin{proof}
In the geometric decomposition of Proposition~\ref{prop:geometric-orbit},
the $e$ summands are conjugate. The idempotent $e_i$ is defined over
$K$, so each summand has the same total dimension over the sheets
above $i$. Their sum is $d_i\beta_i$. This proves the divisibility.
\end{proof}

\begin{remark}\label{rem:kronecker}
A field that splits the algebra need not split an indecomposable
representation. For the Kronecker quiver over $K$, let both vertex
spaces be $K^m$ and take
\[ \quiverdiagram{\xymatrix@C=4em{K^m\ar@/^0.29pc/[r]^{I_m}\ar@/_0.29pc/[r]_{C_m}&K^m}}, \]
where $C_m$ is the companion matrix of $x^m-t$. The endomorphism
algebra is $K_m$.
Replacing $C_m$ by the single Jordan matrix $J_m(c)$, with $c\in K$,
gives residue endomorphism field $K$. Here $I_m$ is the identity
matrix, $C_m$ is chosen so that $K[C_m]\iso K[x]/(x^m-t)$, and
$J_m(c)$ has eigenvalue $c$ and ones on its superdiagonal. Thus
imaginary rank alone does not determine the residue degree.
\end{remark}

\subsection{Ranks of indecomposables and real roots}\label{sec:root}
We now apply the folding results of Section~\ref{sec:unfold}.
The comparison of the dimension period and the isomorphism period
in the proof follows \cite[Proposition~10.2]{DengDu}, with Galois
conjugation in place of Frobenius. Imaginary-root infinitude will follow
from the arithmetic classification in Section~\ref{sec:classification}.

\begin{proposition}\label{prop:root}
The indecomposable ranks of $\Lambda$ are exactly $\Delta^+(C)$.
Every positive real root $\beta$ has a unique indecomposable, whose
residue endomorphism field is $K_{q_d(\beta)}$.
\end{proposition}
\begin{proof}
Let $M$ be indecomposable, and write
$M_{\Kbar}=\bigoplus_{j=1}^e X_j$ as in
Proposition~\ref{prop:geometric-orbit}, and put $\alpha=\dimv X_1$. If $s$ is the least
$a$-period of $\alpha$, equivariance gives $s\mid e$ and
\[ \iota(\rk M)=\frac es\mc N_s^a\alpha. \]
By Kac's theorem and Hubery's Proposition~4,
$\gamma=f(\mc N_s^a\alpha)$ is a positive folded root. If
$\alpha$ is real, uniqueness forces the module and dimension periods
to agree, so $e=s$. If $\alpha$ is imaginary, so is $\gamma$, and its
positive multiple $(e/s)\gamma$ is a root. This proves one inclusion.

Conversely, Hubery's Theorem~1 supplies a constant-field orbit sum
$X=\bigoplus_{j=0}^{s-1}{}^{a^j}Y$ with $\dimv X=\iota\beta$ for every
$\beta\in\Delta^+(C)$. Scalar extension preserves indecomposability:
the local endomorphism algebra has residue field $k$, which becomes
$\Kbar$. It also reflects isomorphisms. Indeed, the determinant of a
general element of the $k$-space $\Hom(Y,Y')$ is a polynomial over
$k$; if it is nonzero after extension, it is nonzero at a $k$-point.
Thus the geometric Galois orbit is exactly the $a$-orbit of $Y$.
Proposition~\ref{prop:descent} descends it to rank $\beta$.

For real $\beta$, Hubery's theorem determines the unique geometric
root orbit and gives its length $q_d(\beta)$. Ordinary real-root
uniqueness and Lemma~\ref{lem:ND} give uniqueness over $K$, and
Proposition~\ref{prop:geometric-orbit} gives the stated residue field.
\end{proof}

\subsection{Euler form}
For rank vectors $\mathbf r,\mathbf s$, put
\[ \langle\mathbf r,\mathbf s\rangle_A =\sum_i d_i r_i s_i-\sum_{i,j}(\dim_K A_{ij})r_i s_j. \]
The usual projective tensor-algebra resolution, as in the ordinary
Hom--Ext calculation of \cite[Chapter~2]{DWbook},
\[ 0\longrightarrow M\otimes_R A\otimes_R\Lambda \longrightarrow M\otimes_R\Lambda\longrightarrow M\longrightarrow0 \]
gives
\[ \dim_K\Hom_\Lambda(M,N)-\dim_K\Ext^1_\Lambda(M,N) =\langle\rk M,\rk N\rangle_A, \qquad \langle\beta,\beta\rangle_A=q_d(\beta). \]
Write $h_M=\dim\End(M)$ for the endomorphism dimension over
the base field (either $K$ or $\mb F_q$). The Euler form is compatible
with the ordinary Euler form $\chi_{\Qbar}$ of the unfolding:
\[
\langle\beta,\gamma\rangle_A
=\chi_{\Qbar}(\iota\beta,\iota\gamma),\qquad
\langle\beta,\gamma\rangle_A+\langle\gamma,\beta\rangle_A
=\beta^{\mathsf T}DC\gamma.
\]
In particular,
\[
q_d(\beta)=\dim G_\beta-\dim V_\beta
=\lVert\beta\rVert_d^2-\dim V_\beta,\qquad
\dim_K\Ext^1_\Lambda(M,M)=h_M-q_d(\beta).
\]
The self-pairing is independent of the orientation; the Euler
form itself need not be symmetric.

\begin{remark}\label{rem:rigidity}
For the real-root indecomposable $M_\beta$, the Euler identity gives
\[ \dim_K\Ext^1_\Lambda(M_\beta,M_\beta)=\dim_K\rad\End_\Lambda(M_\beta). \]
Uniqueness need not imply rigidity. The representation
\[
\quiverdiagram{\xymatrix@C=3.8em@R=2.3em{
 &K^2\ar[dr]^{(0\;1)}&\\
 K\ar[ur]^{\binom10}\ar[rr]_{1}&&K
}}
\]
has rank $(1,2,1)$ on the vertices $1,2,3$, with the middle vertex
at the top. Its endomorphism algebra is $K[z]/(z^2)$ and its
self-extension space has dimension one. This is why
the rigid correspondence of \cite{GLSrigid} uses real Schur roots.
\end{remark}

\subsection{Generic Schur summands}
We call a nonzero representation \emph{Schur} if its endomorphism
algebra is a division algebra, and \emph{absolutely Schur} if its
endomorphism algebra is $K$. By Lemma~\ref{lem:basechange}, the latter
condition is equivalent to having only scalar endomorphisms over
$\Kbar$. The following consequence of Kac's canonical decomposition
and the Happel--Ringel lemma describes general representations of
the ramified species.

\begin{proposition}\label{prop:generic-schur-summands}
For every nonzero rank $\beta$, there is a nonempty Zariski-open subset
$U_\beta\subseteq V_\beta$, defined over $K$, such that every
indecomposable summand $P$ of every $M\in U_\beta(K)$ has
\[ \End_\Lambda(P)\iso K_e \]
for some $e\geq1$. In particular, every indecomposable representation
in $U_\beta(K)$ is Schur.
\end{proposition}
\begin{proof}
After scalar extension to $\Kbar$, Proposition~\ref{prop:unfold}
identifies $V_\beta$ with $\Rep(\Qbar,\iota\beta)$.
Kac's canonical-decomposition theorem
\cite[Section~4 and Proposition~3]{Kac82}, recalled in
\cite[Theorem~2.13]{DWcomb}, gives a dense open subset on which the
indecomposable summands are Schur and the extension groups between
nonisomorphic indecomposable summands vanish. Equal dimension vectors
are allowed among these summands.
This open subset is defined over a finite Galois extension of $K$.
Intersecting its finitely many Galois conjugates gives a nonempty
open subset that descends to $U_\beta$.

Let $P$ be an indecomposable summand of $M\in U_\beta(K)$.
Proposition~\ref{prop:geometric-orbit} gives
\[ P_{\Kbar}\iso X_1\oplus\cdots\oplus X_e, \]
where the $X_j$ are pairwise nonisomorphic and form one Galois orbit.
By Krull--Schmidt, they occur among the geometric indecomposable
summands of $M$. Thus each $X_j$ is Schur, and
$\Ext^1(X_l,X_j)=0$ for $j\ne l$.
By the Happel--Ringel lemma \cite[Lemma~4.1]{HR82}, every nonzero map
$X_j\to X_l$ is a monomorphism or an epimorphism. The $X_j$ have the
same total $\Kbar$-dimension, since Galois conjugation permutes the
unfolded vertices. Such a map would therefore be an isomorphism.
Consequently $\Hom(X_j,X_l)=0$ for $j\ne l$, and
\[ \Kbar\otimes_K\End_\Lambda(P)
   \iso\End_{\Lambda_{\Kbar}}(P_{\Kbar})\iso\Kbar^{\,e}. \]
Lemma~\ref{lem:basechange} implies that $\End_\Lambda(P)$ has zero
radical. It is local, hence a division algebra, and
Lemma~\ref{lem:fields} identifies it with $K_e$.
\end{proof}

For an indecomposable $M\in U_\beta(K)$, its residue degree $e$
therefore satisfies
\[ e\mid g_d(\beta),\qquad
   \dim_K\Ext^1_\Lambda(M,M)=e-q_d(\beta). \]
The first assertion is Proposition~\ref{prop:divisibility}; the
second follows from the Euler identity.

\begin{corollary}\label{cor:absolute-schur-ranks}
For a nonzero rank $\beta$, the following are equivalent:
\begin{enumerate}
\item $\iota\beta$ is a Schur root of $\Qbar$.
\item A general geometric representation in $V_\beta$ is indecomposable.
\item There is $M\in V_\beta(K)$ with $\End_\Lambda(M)=K$.
\end{enumerate}
When they hold, the absolutely Schur locus is a nonempty
Zariski-open subset of $V_\beta$, defined over $K$.
\end{corollary}
\begin{proof}
A Schur root is a dimension vector admitting a representation with
only scalar endomorphisms. The first two conditions are equivalent
by \cite[Proposition~1(a)]{Kac82}; see also
\cite[Section~2.3]{DWcomb}. The locus of endomorphism dimension
one is open and defined over $K$. If it is geometrically nonempty,
it has a $K$-point, since $V_\beta$ is an affine space and $K$ is
infinite. Conversely, a representation with endomorphism algebra $K$
has only scalar endomorphisms after base change to $\Kbar$.
\end{proof}

\begin{remark}\label{rem:arithmetic-schur-genericity}
Proposition~\ref{prop:generic-schur-summands} concerns a dense open
subset of the whole representation space. The word \emph{absolute} in
Corollary~\ref{cor:absolute-schur-ranks} is also essential. In
Remark~\ref{rem:kronecker}, $(I_2,C_2)$ has endomorphism field $K_2$,
but the geometric canonical decomposition of $(2,2)$ is
$(1,1)\oplus(1,1)$. Moreover, decomposable $K$-representations of this
rank are Zariski dense: on the open set where the first arrow is
invertible, their arrow ratio can be any split regular semisimple
matrix. Conjugates of diagonal matrices with distinct eigenvalues
are dense over the infinite field $K$. Thus existence of a Schur
representation over $K$ does not imply that any nonempty open
subset has only indecomposable $K$-points.
\end{remark}

\section{Deformed preprojective algebras and generic parameters}
\label{sec:moment}

Fix a rank vector $\beta$ and the vertex spaces
$M_i=K_{d_i}^{\beta_i}$. Representations on these fixed vertex spaces form the affine space
\[ V_\beta=\bigoplus_{i,j} \Hom_{K_{d_j}}(M_i\otimes_{K_{d_i}}A_{ij},M_j). \]
The group $G_\beta=\prod_i\Res_{K_{d_i}/K}\mr{GL}_{\beta_i}$
acts by change of vertex bases. Here Weil restriction describes this
as an algebraic group over $K$; in particular,
$G_\beta(K)=\prod_i\mr{GL}_{\beta_i}(K_{d_i})$.
Its Lie algebra is
$\mathfrak g_\beta=\bigoplus_i\End_{K_{d_i}}(M_i)$, regarded as a
$K$-vector space. This is distinct from the Kac--Moody root space
$\mathfrak g(C)_\beta$. Writing $M_{ij}$ for an arrow map, the
differential $d_M:\mathfrak g_\beta\to V_\beta$ is
\[ (d_Mf)_{ij}=f_jM_{ij}-M_{ij}(f_i\otimes1). \]
We also denote this infinitesimal action by $f\cdot M$.
We use the nondegenerate \emph{field-trace pairing}
\[ (f,h)=\sum_i\Tr_{K_{d_i}/K} \bigl(\op{tr}_{K_{d_i}}(f_i h_i)\bigr) \quad(f,h\in\mathfrak g_\beta). \]
Here $\op{tr}_{K_{d_i}}$ is the ordinary trace of a
$K_{d_i}$-linear endomorphism, and $\Tr_{K_{d_i}/K}$ is the field
trace from $K_{d_i}$ to $K$.
The multiplicative group $\mb G_m$ acts by the same scalar at
every vertex, so it acts trivially on $V_\beta$. A star denotes the
$K$-linear dual. We write
$\mathfrak g_\beta^{*,0}=\{\ell\in\mathfrak g_\beta^*:\ell(1_M)=0\}$,
where $1_M=(1_{M_i})_i$ is the scalar identity.
The \emph{cotangent moment map} is
\[ \mu_\beta:V_\beta\oplus V_\beta^*\longrightarrow\mathfrak g_\beta^{*,0}, \qquad \mu(M,\xi)(f)=\xi(f\cdot M). \]
When $\beta$ is fixed, we abbreviate $\mu_\beta$ to $\mu$ and write
$\mu_M(\xi)=\mu(M,\xi)$.
For $\boldsymbol\lambda=(\lambda_i)\in\prod_iK_{d_i}$, its
parameter functional on $\mathfrak g_\beta$ is
\[ \ell_{\boldsymbol\lambda}(f)= \sum_i\Tr_{K_{d_i}/K} \bigl(\op{tr}_{K_{d_i}}(\lambda_i f_i)\bigr). \]
Throughout, $\mu=\boldsymbol\lambda$ abbreviates
$\mu=\ell_{\boldsymbol\lambda}$.
It has no solution unless this functional annihilates the global
scalar identity.

After extension to $K_N$, these are the ordinary cotangent
representation space and moment map for $\Qbar$. The parameters at
the sheets above $i$ are the conjugates of $\lambda_i$.
Let $\Qbar^{\mr{dbl}}$ be the \emph{double quiver}, obtained by
adjoining a reverse arrow $b^*$ for each arrow
$b:t(b)\to h(b)$ of $\Qbar$:
\[ \quiverdiagram{\xymatrix@C=3.8em{i\ar[r]^{b}&j}}\qquad\longmapsto\qquad
\quiverdiagram{\xymatrix@C=3.8em{i\ar@/^0.6pc/[r]^{b}&j\ar@/^0.6pc/[l]^{b^*}}}. \]

With our left-to-right path convention,
the moment relations are
\[ \rho_v=\sum_{h(b)=v}b^*b-\sum_{t(b)=v}bb^*. \]
Indeed, the coefficient of $f_v$ in $\xi(d_Mf)$ is the trace
pairing with the operator represented by $\rho_v$. If
$\widetilde\lambda_v$ is the conjugate parameter at the sheet $v$,
set
\[ \Pi_A^{\boldsymbol\lambda} =\left(K_N\Qbar^{\mr{dbl}}/ (\rho_v-\widetilde\lambda_v e_v)_v\right)^{\Gal(K_N/K)}. \]
Here Galois conjugation acts on both the coefficients and the quiver, as in
Proposition~\ref{prop:unfold}. It preserves the relations, so Galois
descent identifies its scalar extension to $K_N$ with the quotient above.
Representations of $\Pi_A^{\boldsymbol\lambda}$ on the fixed vertex
spaces are exactly the solutions of $\mu=\ell_{\boldsymbol\lambda}$.
We call $\Pi_A^{\boldsymbol\lambda}$ the
\emph{deformed preprojective algebra} of the species; for
$\boldsymbol\lambda=0$ we write $\Pi_A$.
These preprojective algebras belong to the generic hereditary species
and must not be confused with the GLS preprojective algebras on the
special fibre.
Reversing an arrow replaces its bimodule by the trace-pairing dual;
a sign change on one arrow of the doubled pair identifies the
corresponding moment-map equations.

\begin{proposition}\label{prop:moment}
For $M\in V_\beta$, the fibre $\mu_M^{-1}(\boldsymbol\lambda)$ is
nonempty precisely when $\ell_{\boldsymbol\lambda}$ annihilates
$\End_\Lambda(M)$. If nonempty, it is an affine space over $K$ of
dimension $\dim_K\Ext^1_\Lambda(M,M)=h_M-q_d(\beta)$.
For parameters $\lambda\in K^I$, this criterion is equivalent to
\begin{equation}\label{eq:weightedcriterion} \lambda^{\mathsf T}D\rk N=0 \quad\text{for every indecomposable direct summand }N\text{ of }M.\end{equation}
\end{proposition}
\begin{proof}
We apply the standard resolution to obtain
\[ 0\longrightarrow\End_\Lambda(M)\longrightarrow\mathfrak g_\beta \xrightarrow{d_M}V_\beta\longrightarrow\Ext^1_\Lambda(M,M) \longrightarrow0. \]
Dualizing over $K$, and viewing $\mu_M$ as a map to
$\mathfrak g_\beta^*$, gives
\[ 0\longrightarrow\Ext^1_\Lambda(M,M)^*\longrightarrow V_\beta^* \xrightarrow{\mu_M}\mathfrak g_\beta^* \longrightarrow\End_\Lambda(M)^*\longrightarrow0. \]
This proves the first assertion and the dimension formula.

Suppose now that $\lambda_i\in K$. Necessity of
\eqref{eq:weightedcriterion} follows by applying the functional to the
idempotent projecting onto $N$. For sufficiency, extend to $\Kbar$.
An indecomposable $N$ becomes one multiplicity-free Galois orbit
$X_1\oplus\cdots\oplus X_e$. Write $\bar\lambda$ for the parameter
constant with value $\lambda_i$ on the sheets over $i$. The values
$\bar\lambda\cdot\dimv X_j$ are equal: conjugation only permutes
sheets within the original vertices, and the parameters are fixed.
Their sum is $\lambda^{\mathsf T}D\rk N$. Since the characteristic is
zero, each value vanishes if and only if the sum vanishes.
The ordinary criterion \cite[Theorem~8.1.3]{DWbook} therefore gives a
solution over $\Kbar$. Equivalently, the functional annihilates the
geometric endomorphism algebra. This is a linear system defined over
$K$, so it has a solution over $K$. Applying this to the summands of
$M$ proves sufficiency.
\end{proof}

The first part of Proposition~\ref{prop:moment} remains the correct
criterion for general $\lambda_i\in K_{d_i}$. It cannot be replaced
by the vanishing of the total trace alone. At a single vertex with field
$K_2=K(u)$, $u^2=t$, and no arrows, take $M=K_2$ and parameter $u$.
Its split quiver consists of two isolated vertices $v_0$ and $v_1$
with $(\widetilde\lambda_{v_0},\widetilde\lambda_{v_1})=(u,-u)$.
Although $\Tr_{K_2/K}(u)=0$, the moment map vanishes identically
while the parameter functional does not. Hence the moment-map equation
has no solution.

\subsection{Invariant generic parameters}

A parameter for the split quiver is a vector
$\theta\in\Kbar^{\Qbar_0}$. We use the ordinary coordinate pairing
$\theta\cdot\gamma=\sum_v\theta_v\gamma_v$. For
$\alpha\in\mb N^{\Qbar_0}$, the parameter is
\emph{$\alpha$-generic} if $\theta\cdot\alpha=0$ and
$\theta\cdot\gamma\ne0$ for every integral vector
$0<\gamma<\alpha$. The inequalities here are coordinatewise; some
coordinates of a proper subvector may equal those of $\alpha$.
A parameter constant on vertex orbits has the form
$\bar\lambda=\iota\lambda$.

\begin{proposition}\label{prop:genericity}
Let $\beta\ne0$ and $\alpha=\iota\beta$. There exists a parameter
$\lambda\in\mb Z^I$ for which $\bar\lambda$ is
$\alpha$-generic if and only if $g_d(\beta)=1$.
\end{proposition}
\begin{proof}
We set $v_i=d_i\beta_i$. The invariant parameter space satisfying the
total equation is the rational hyperplane $\lambda\cdot v=0$.
For a subvector $\gamma\leq\alpha$, put
$n_i=\sum_{r\in\mb Z/d_i\mb Z}\gamma_{(i,r)}$.
Its equation on this hyperplane is $\lambda\cdot n=0$.
This equation holds identically if and only if $n$ is a rational
multiple of $v$.

If $g_d(\beta)=1$, the integral vector $v$ is primitive. The only integral vectors in its rational span lying coordinatewise
between $0$ and $v$ are $0$ and $v$. These totals force $\gamma=0$ and $\gamma=\alpha$,
respectively. Hence every proper subvector cuts out a proper
hyperplane. There are finitely many, so a rational parameter avoids
them all. Clear denominators.

Conversely, if $g_d(\beta)>1$, distribute $v_i/g_d(\beta)$ units among the
$d_i$ sheets above $i$, with at most $\beta_i$ units at each sheet.
This is possible since $0\leq v_i/g_d(\beta)\leq d_i\beta_i$. The resulting
nonzero proper integral subvector $\gamma$ has totals $v/g_d(\beta)$.
Every invariant parameter annihilating $\alpha$ also annihilates
$\gamma$, and none is generic.
\end{proof}

In particular, primitive rank need not suffice: a vertex simple
with $d_i>1$ admits no invariant generic parameter.

\subsection{Generic parameters in the full vertex fields}
\label{subsec:full-parameters}

We now allow the parameter at $i$ to lie in $K_{d_i}$ rather than in
$K$. After scalar extension, its coordinates above $i$ are its field
conjugates. The resulting vector is fixed by the semilinear Galois
action, although it need not be fixed by the quiver automorphism alone.
This allows more parameters than the invariant construction of
Proposition~\ref{prop:genericity}.

For $\theta=(\theta_i)\in R$, let $\widetilde\theta$ be its vector of
split coordinates. Thus
$\widetilde\theta_{(i,s)}=\sigma_{i,s}(\theta_i)$ for the $K$-embeddings
$\sigma_{i,s}:K_{d_i}\hookrightarrow\Kbar$. The \emph{trace hyperplane} is
\[ \mc H_\beta=\left\{\theta\in R: \sum_i\beta_i\Tr_{K_{d_i}/K}(\theta_i)=0\right\}. \]
This is a $K$-linear hyperplane when $\beta\ne0$. Its scalar extension
is the full hyperplane $\widetilde\theta\cdot\iota\beta=0$.

\begin{proposition}\label{prop:full-generic}
For $\beta\ne0$, there is a $\theta\in \mc H_\beta$ for which
$\widetilde\theta$ is generic for $\iota\beta$ if and only if
$\gcd_i\beta_i=1$. When this condition holds, $\theta$ can be chosen
in $\prod_i k[[t^{1/d_i}]]$.
\end{proposition}
\begin{proof}
After scalar extension, $\mc H_\beta$ is the hyperplane orthogonal
to $\iota\beta$. The linear form associated with a subvector vanishes identically on
this hyperplane precisely when the subvector is proportional to
$\iota\beta$. Such a proper nonzero integral subvector exists
precisely when $\beta$ is divisible.
Otherwise the hyperplanes defined by proper subvectors, which we
call the \emph{walls}, form a finite Galois-stable union of proper
hyperplanes. Its complement is a nonempty open subset of the
$K$-affine space $\mc H_\beta$, and hence has a $K$-point.
Multiplying by a sufficiently large power of $t$ makes every vertex
coordinate integral without changing which pairings vanish.
\end{proof}

For example, orient the degree-$(1,4)$ species from the vertex with
field $K$ to the vertex with field $K_4$. Its species diagram and
unfolding are
\begin{equation} \label{eq:1,4}
\quiverdiagram{\xymatrix@C=3em{K\ar[r]^{K_4}&K_4}}
\qquad\qquad
\quiverdiagram{\xymatrix@C=2.2em@R=1.7em{
 &v_0&\\
 v_3&c\ar[u]\ar[r]\ar[d]\ar[l]&v_1\\
 &v_2&
}},\qquad a=(v_0\ v_1\ v_2\ v_3).
\end{equation}

The automorphism fixes $c$. At rank $(2,1)$ the dimensions are two
at $c$ and one at each leaf. Write $u^4=t$. The parameter
$\theta_1=-2,\ \theta_2=1+u+u^2+u^3$
is generic for the unfolded vector $(2,1,1,1,1)$. To check this,
let $n\in\{0,1,2\}$ and $n_s\in\{0,1\}$, $0\le s<4$, be a
subvector annihilated by the parameter. If $\zeta$ is a primitive
fourth root of unity, the coefficients of $u,u^2,u^3$ give
$\sum_{s=0}^3 n_s\zeta^{sr}=0$ for $r=1,2,3$.
The Fourier inversion formula forces the four $n_s$ to be equal.
The constant coefficient then gives $2n=\sum_s n_s$, so the
subvector is zero or the full vector. Here $g_d(2,1)=2$; the
parameter is preserved by the combined action on coefficients and
vertices, although it is not fixed by the vertex permutation alone.
Corollary~\ref{cor:full-generic-moduli} describes the corresponding
moment-map quotient.

\subsection{Reflections}

\begin{proposition}\label{prop:reflections}
At a sink or source $i$, the species reflection functor gives an
equivalence between the full subcategories having no direct summand
isomorphic to the simple at $i$. It commutes with field extension,
preserves endomorphism algebras, and sends rank $\beta$ to $s_i\beta$.
\end{proposition}
\begin{proof}
The vertices above $i$ are pairwise nonadjacent. The ordinary
reflection functors at these vertices commute. Their product
commutes with the semilinear Galois action: the action permutes
the vertices, and the kernel and cokernel constructions commute
with scalar extension. It therefore descends, together with its
inverse, to the stated equivalence of species representations.
The product of the simple reflections above $i$ restricts to $s_i$
on invariant dimension vectors. The remaining assertions follow
from the equivalence and its construction.
\end{proof}

\section{Absolutely indecomposable families}\label{sec:imaginary}
\label{sec:absolute-families}
We now construct the families needed to prescribe the Galois orbit
of an indecomposable representation. A variety is \emph{geometrically
integral} if it remains reduced and irreducible over an algebraic
closure. We work with smooth geometrically integral,
$G_\beta$-invariant locally closed subvarieties
$W\subseteq V_\beta$ that have a $K$-point, have constant endomorphism
dimension $h$, and consist of geometrically indecomposable
representations. We call such a $W$ a \emph{family of absolutely
indecomposable representations}. In this section, ``family'' always
includes these conditions.

Every orbit in $W$ has dimension $\dim G_\beta-h$, so
\[ \nu(W)=\dim W-\lVert\beta\rVert_d^2+h \]
is the dimension of a rational quotient of $W$ by $G_\beta$. Equivalently,
Rosenlicht's theorem gives
$\nu(W)=\op{trdeg}_K K(W)^{G_\beta}$, where $K(W)$ is the
field of rational functions and the superscript denotes invariants.
This is the number of
parameters after change of bases, rather than the dimension of the
representation space itself.

For an ordinary quiver, write $p(\alpha)=1-q_Q(\alpha)$. Its Kac
polynomial $A_{Q,\alpha}(T)$ is monic of degree $p(\alpha)$ when
$\alpha$ is a positive root, and is zero otherwise \cite{Kac83,HLRV}.
We use this fact to bound the dimensions of families. Roots with
$q_Q(\alpha)<0$ and roots with $q_Q(\alpha)=0$ require different
constructions.

\subsection{Rational points and reflections}
For a finite extension $F=k((s))$ of $K$, we use the valuation topology:
a basic neighbourhood of $x\in F$ is $x+s^n k[[s]]$. On rational
points of a variety it induces the corresponding topology in local
affine coordinates. We call these neighbourhoods valuation-open.

\begin{lemma}\label{lem:henseldense}
Let $F$ be a finite extension of $K$, and let $Z$ be a smooth
geometrically integral $F$-variety with an $F$-point $z$. Every
sufficiently small valuation-open neighbourhood of $z$ is Zariski
dense in $Z$. If a regular function $\varphi$ is smooth at $z$, its values
on such a neighbourhood contain a valuation-open neighbourhood of
$\varphi(z)$ in $F$.
\end{lemma}
\begin{proof}
An \'etale coordinate map at $z$, chosen to include $\varphi$ when appropriate,
identifies a small neighbourhood with an open polydisc by Hensel's
lemma. A polynomial vanishing on an open polydisc is zero, by induction
on the number of variables. A proper closed subset of $Z$ has smaller
dimension, and its image under an \'etale map cannot contain that
polydisc. The last assertion follows by projection to the coordinate
$\varphi$.
\end{proof}

\begin{lemma}\label{lem:family-reflection}
Suppose $\iota\beta$ is an imaginary root and $i$ is a sink or source.
Reflection sends a family of rank $\beta$ to a family of rank
$s_i\beta$ with the same value of $\nu$.
\end{lemma}
\begin{proof}
We treat a sink; the source case is dual. At every sheet over the
sink, the incoming map is surjective: otherwise a sink simple would
split off from a geometric indecomposable. The kernels have constant
dimension. Their \emph{frames}, meaning choices of ordered bases,
form a principal bundle over $W$. A frame identifies each kernel with the fixed vertex
space of the reflected dimension and hence defines a point of the
reflected representation space. Conversely, the outgoing map in a
reflected representation is injective. Choosing a frame of its
cokernel reconstructs the original representation. The two constructions commute with field extension.

The bundle of kernel frames has fibre dimension $d_i(s_i\beta)_i^2$.
Over its image, the reverse frame choices form torsors of dimension
$d_i\beta_i^2$ under $\Res_{K_{d_i}/K}\mr{GL}_{\beta_i}$;
these torsors split over $K$ by Hilbert~90. We take a smooth invariant
open subset of the image closure contained in the image. Its inverse
image in the frame bundle has a rational point by
Lemma~\ref{lem:henseldense}. Reflection preserves endomorphisms and
geometric indecomposability, so $W'$ is a family in the sense fixed
at the start of this section. The frame dimensions give
\[ \dim W'-\dim W=d_i\bigl((s_i\beta)_i^2-\beta_i^2\bigr) =\dim G_{s_i\beta}-\dim G_\beta. \]
Since the endomorphism dimension is unchanged, $\nu(W')=\nu(W)$.
\end{proof}

\subsection{Two consequences of the ordinary counting theorem}

We apply the Kac-polynomial degree bound above to orbit quotients.
In the isotropic case, we also obtain a rational curve of
isomorphism classes. The details below ensure compatibility with
descent and quiver automorphisms.

\begin{lemma}[Orbit quotients and parameter bounds]
\label{lem:orbit-parameter-bound}
Let $Q$ be a finite acyclic quiver, let $\alpha$ be a positive root,
and put $p(\alpha)=1-q_Q(\alpha)$. Over an algebraically closed field
of characteristic zero, the indecomposable locus of dimension $\alpha$
has a finite stratification by invariant locally closed varieties
$S_i$ with \emph{geometric quotients} $\pi_i:S_i\to Y_i$ for the
change-of-basis action. In particular, their geometric fibres are single
isomorphism orbits, and invariant morphisms factor through $\pi_i$. We may choose $S_i$ and $Y_i$
smooth and geometrically integral, with constant endomorphism dimension
on each $S_i$. Every $Y_i$ has dimension at most $p(\alpha)$.
Consequently, an invariant family with constant endomorphism dimension
$h$ has dimension at most
\[ \dim G_\alpha-h+p(\alpha). \]
If $p(\alpha)=1$, one quotient is an open subset of a smooth projective
rational curve and its stratum contains all but finitely many
indecomposable orbits. After removing finitely many further orbits,
this quotient can be made equivariant for any finite group of quiver
automorphisms fixing $\alpha$.
\end{lemma}
\begin{proof}
We first work over $\overline{\mb Q}$. The indecomposable locus
is constructible: its complement is the union of the images of the
finitely many direct-sum morphisms. Rosenlicht's theorem
\cite{Rosenlicht}, applied successively to the complement of each open
stratum, gives finitely many geometric quotients $\pi_i:S_i\to Y_i$.
We refine the strata by irreducible components, smooth loci, and
endomorphism dimension. These refinements are invariant under the
connected group $G_\alpha$. Applying Rosenlicht again to their
complements gives the asserted regularity properties.

All this involves finitely many varieties and morphisms defined over
one number field. We spread them over a localization of its ring of
integers. After deleting finitely many primes, the strata still
partition the geometric indecomposable locus, the maps $\pi_i$ are
surjective, and their geometric fibres are single orbits. For the last
assertion it suffices to preserve the surjectivity of
\[ G_\alpha\times S_i\longrightarrow S_i\times_{Y_i}S_i,\qquad (g,x)\longmapsto(gx,x). \]
The complement of its image is constructible and has empty geometric
generic fibre. Its image on the one-dimensional arithmetic base is
therefore finite. The same argument preserves surjectivity of $\pi_i$
and the finite partition, using the relative direct-sum morphisms.
Smoothness, geometric integrality, and the relevant dimensions hold
after further shrinking. We use these properties of the chosen
quotient morphisms in the counting argument below.

Fix a remaining residue field $\mb F_{q_0}$. Over every finite
extension of it, a rational point of $Y_i$ specifies a Frobenius-stable
geometric indecomposable class. Finite-field descent gives a rational
representative. It is unique up to rational isomorphism, since the
stabilizer is the unit group of an endomorphism algebra, hence a
geometrically connected group. This is also the fixed-orbit case of
\cite[Lemma~3.1, Theorem~5.1, and Corollary~5.4]{DengDu}. Thus
\begin{equation}\label{eq:quotient-count-spread} A_{Q,\alpha}(q_0^n)=\sum_i |Y_i(\mb F_{q_0^n})|\qquad(n\ge1).\end{equation}
Each summand has leading term $q_0^{n\dim Y_i}$ by Lang--Weil
\cite{LangWeil}. The Kac polynomial is monic of degree $p(\alpha)$.
Consequently exactly one quotient has dimension $p(\alpha)$, and none
has larger dimension. Intersecting any invariant family with the
strata and using orbit dimension $\dim G_\alpha-h$ proves its stated
dimension bound. Base extension gives the same conclusions over every
algebraically closed field of characteristic zero.

Suppose $p(\alpha)=1$. There is one curve quotient $U$, while all
other quotients consist of finitely many points. Let $\overline U$ be
its smooth projective completion. Spread this completion as a smooth
proper curve as well. After extending the chosen finite field, all
boundary points and all zero-dimensional quotients are rational.
Equation~\eqref{eq:quotient-count-spread}, with $A_{Q,\alpha}(T)=T+c$,
then gives
\[ |\overline U(\mb F_{q_0^n})|=q_0^n+b\qquad(n\ge1) \]
for a fixed integer $b$. Comparing zeta functions with the Weil formula
\cite[Sections~14 and 26--27]{MilneEtale} yields
\[ \prod_{j=1}^{2g}(1-\lambda_j z)=(1-z)^{1-b},\qquad |\lambda_j|=q_0^{1/2}>1. \]
The left side is a polynomial with no zero at $z=1$. Thus $b=1$ and
$g=0$. Genus is constant in the smooth proper model, so the
characteristic-zero completion is a rational curve.

Finally let a finite group act on $Q$ and fix $\alpha$. The stratum
over $U$ contains all but finitely many indecomposable orbits. Its
intersection with its finitely many translates is invariant and still
has this property. Each orbit is a fibre of the geometric quotient,
so deleting these finitely many orbits amounts to deleting finitely
many points of $U$. On this common domain the quotient property
induces a birational action on $U$. The group identities hold on a
dense open subset, and each birational map extends uniquely to an
automorphism of the smooth projective completion. Delete the finite
orbits of all excluded points. The resulting open curve and its
inverse image give an equivariant geometric quotient. Further
invariant open shrinkings make the source smooth with constant
endomorphism dimension, without losing the cofinite property.
\end{proof}

The quotient in Lemma~\ref{lem:orbit-parameter-bound} parametrizes
orbits on a locally closed part of the indecomposable locus. It is not
the affine invariant-theoretic quotient of $V_\alpha$, which is a point
when the quiver is acyclic. The geometric quotients are useful here
because their points correspond to individual isomorphism classes.

\subsection{The case $q_d(\beta)<0$: changing orientation}

Write $\beta=m\beta_0$ with $\beta_0$ primitive, and put
$\alpha=\iota\beta=m\alpha_0$. At divisible dimension, every parameter
orthogonal to $\alpha$ is also orthogonal to $\ell\alpha_0$.
We therefore choose a parameter whose only vanishing pairings are the
ones forced by this proportionality:
\begin{equation}\label{eq:proportional-parameter} \widetilde\theta\cdot\eta=0 \quad\Longleftrightarrow\quad \eta=\ell\alpha_0\text{ for some }0\le\ell\le m, \qquad 0\le\eta\le\alpha,\quad \eta\in\mb Z^{\Qbar_0}.\end{equation}
The next lemma shows that parameters in the full vertex fields allow
us to avoid every other subdimension vector.

\begin{lemma}\label{lem:ray-generic}
For every nonzero $\beta$, there is a parameter
$\theta\in\mc H_\beta$ satisfying
\eqref{eq:proportional-parameter}.
\end{lemma}
\begin{proof}
We apply the hyperplane argument of Proposition~\ref{prop:full-generic}
to the primitive vector $\alpha_0$, considering all integral subvectors
of $\alpha=m\alpha_0$. On the split trace hyperplane, precisely the
proportional subvector equations vanish identically. Every other one
cuts out a proper hyperplane. The complement of their finite,
Galois-stable union is a nonempty $K$-open subset of $\mc H_\beta$,
so it contains a $K$-point. The proportional integral subvectors are
exactly $\ell\alpha_0$, $0\le\ell\le m$.
\end{proof}

\begin{lemma}\label{lem:negative-transport}
Suppose $\alpha=\iota\beta$ is an imaginary root with
$q_{\Qbar}(\alpha)<0$. If one acyclic orientation admits a family $W$ of rank $\beta$
with $\nu(W)=p(\alpha)$, then every
acyclic orientation of the same modulation has such a family.
\end{lemma}
\begin{proof}
Write $\alpha=m\alpha_0$, set $r=-q_{\Qbar}(\alpha_0)\ge1$,
and choose $\theta$ as in Lemma~\ref{lem:ray-generic}. Its trace functional vanishes
on the identity and on the nilpotent endomorphisms of every geometric
member of $W$, hence on its whole endomorphism algebra. By
Proposition~\ref{prop:moment}, the space of choices of reverse arrows
$Z_W\to W$ is an affine bundle of rank $h-q_{\Qbar}(\alpha)$.
Since $q_{\Qbar}(\alpha)=q_d(\beta)
=\dim G_\beta-\dim V_\beta$, the endomorphism and
change-of-bases dimensions cancel:
\begin{equation}\label{eq:negative-lift-dimension}
\dim Z_W=\dim V_\beta+p(\alpha).
\end{equation}
This bundle is smooth, geometrically integral, and has a $K$-point.

We project to the hereditary representation space for the new
orientation. If a restriction decomposes into $\kappa\ge2$ indecomposables,
the lifting criterion and \eqref{eq:proportional-parameter} force their
dimensions to be
$\ell_1\alpha_0,\ldots,\ell_\kappa\alpha_0$, with
$\sum_j\ell_j=m$. Choose geometric quotient strata for these summands as in
Lemma~\ref{lem:orbit-parameter-bound}. Let their quotient dimensions be
$\delta_j\le p(\ell_j\alpha_0)$. A direct sum is determined by the
isomorphism classes of its summands, up to a finite permutation. Thus
these sums depend on at most $\sum_j\delta_j$ parameters. On the
stratum where the total endomorphism dimension is $h'$, each orbit has
dimension $\dim G_\beta-h'$, and the locus of such sums has dimension at most
\[ \dim G_\beta-h'+\sum_j\delta_j. \]
The choices of reverse arrows over such a representation form an
affine space of dimension $h'-q_{\Qbar}(\alpha)$. Adding the two
dimensions cancels $h'$. Hence the locus of lifts with this
decomposition type has dimension at most
\[
\dim V_\beta+\sum_jp(\ell_j\alpha_0).
\]

The difference from the dimension in
\eqref{eq:negative-lift-dimension} is positive. Indeed,
$m^2-\sum_j\ell_j^2=2\sum_{i<j}\ell_i\ell_j\ge \kappa(\kappa-1)$, so
\begin{align*}
p(m\alpha_0)-\sum_jp(\ell_j\alpha_0)
 &=1-\kappa+r\Big(m^2-\sum_j\ell_j^2\Big)\\
 &\ge(\kappa-1)(r\kappa-1)>0.
\end{align*}
There are only finitely many partitions and endomorphism strata.
Removing the closures of these loci leaves a nonempty open subset
of $Z_W$ whose restrictions are geometrically indecomposable. When
$m=1$, no proper decomposition satisfies \eqref{eq:proportional-parameter}.

Choose a smooth invariant open subset $W'$ of the image closure,
contained in this image, with constant endomorphism dimension $h'$.
Its inverse image is nonempty open in $Z_W$ and has a $K$-point by
Lemma~\ref{lem:henseldense}. The fibre bound gives
$\dim W'\ge\dim W+h-h'$, so $\nu(W')\ge p(\alpha)$.
Lemma~\ref{lem:orbit-parameter-bound} gives the opposite inequality.
Hence $W'$ is the required family, with $\nu(W')=p(\alpha)$.
\end{proof}

When $q_{\Qbar}(\alpha)=0$, the same difference is $1-\kappa$, so this
argument gives no dimension bound excluding decompositions. We treat
that case using the rational curve from
Lemma~\ref{lem:orbit-parameter-bound}.

\subsection{Isotropic roots and a descended rational curve}

The root-theoretic description of isotropic roots involves affine
subdiagrams. Indeed, if a connected vector $\eta$ lies in the
fundamental region and $q_Q(\eta)=0$, then
\[ 0=(\eta,\eta)=\sum_{v\in\supp\eta}\eta_v(\eta,\epsilon_v). \]
Every term is nonpositive, so the Cartan matrix on the support has
the positive null vector $\eta$. That support is affine, and $\eta$
is a multiple of its primitive null root \cite[Chapter~4]{KacBook}.
For a prescribed orientation, however, the Weyl element reaching this
vector need not be realized by sink and source reflection functors.
We use the orbit quotient in Lemma~\ref{lem:orbit-parameter-bound}
to construct the family without changing orientation.

\begin{proposition}\label{prop:isotropic-absolute}
If $\alpha=\iota\beta$ is a positive isotropic root of the
unfolding, with no primitivity assumption, the ramified species
admits a family of rank $\beta$ with $\nu=1$.
\end{proposition}
\begin{proof}
Apply Lemma~\ref{lem:orbit-parameter-bound} to $\alpha$ and the group
$\langle a\rangle$. We obtain an equivariant geometric quotient
$\pi:W\to U$ over $k$, where $U$ is an open subset of a smooth
projective rational curve. The variety $W$ is smooth and geometrically
integral, its endomorphism dimension is constant, and all its geometric
points are indecomposable representations of dimension $\alpha$.

Over $K_N$, combine coefficient conjugation with the action of $a$.
Finite Galois descent gives a quotient $\pi_0:W_0\to U_0$ for the
ramified species. The projective completion of $U_0$ is a $K$-form of
$\mb P^1$. Its class lies in $H^1(K,\mr{PGL}_2)$, which injects
into $\Br(K)$. Lemma~\ref{lem:fields} therefore identifies this
completion with $\mb P^1_K$. In particular $U_0(K)$ is infinite.

Choose $u\in U_0(K)$. The geometric fibre over $u$ is a single orbit
of the change-of-basis group. Because the fibre is defined over $K$,
this isomorphism class is fixed by $\Gamma$. Proposition~\ref{prop:descent}
applied to the one-element Galois orbit gives a $\Lambda$-module $M$
whose scalar extension belongs to that class and is indecomposable.
After choosing bases in the $K_{d_i}$-spaces $Me_i$, we obtain a point
of the representation space over $K$. It lies in $W_0$ and maps to $u$:
both assertions can be checked after extension to $\overline K$, where
they follow from invariance under change of bases and the description
of the fibre. Thus the rational quotient point does give an actual
absolutely indecomposable representation over $K$.

The descended variety is smooth and geometrically integral, with
constant endomorphism dimension $h$. Its fibres have dimension
$\dim G_\beta-h$. Since $\dim U_0=1$, we obtain
\[ \nu(W_0)=\dim W_0-\lVert\beta\rVert_d^2+h=1. \]
The argument also applies over every finite extension of $K$, since
Lemma~\ref{lem:fields} holds over those fields.
\end{proof}

\subsection{The fundamental-region estimate}
The following dimension estimate distinguishes primitive vectors
from those with $q_{\Qbar}(\eta)<0$. It supplies the strict
inequality used to construct families by changing orientation.

\begin{lemma}\label{lem:fundamental-open}
Let $\eta>0$ have connected support and satisfy
$(\widetilde C\eta)_v\le0$ for every vertex. If $\eta$ is primitive
or $q_{\Qbar}(\eta)<0$, the representation space at $\eta$ has a
nonempty open subset consisting of geometrically indecomposable
representations.
\end{lemma}
\begin{proof}
We restrict to the support. For a proper nonzero integral subvector
$\gamma$ put $r_v=\gamma_v/\eta_v$. If $m_{vw}$ is the number of
unoriented edges between $v$ and $w$, expansion of the symmetric
Cartan form gives
\begin{equation}\label{eq:fundamental-identity} (\gamma,\eta-\gamma) =\sum_v\eta_vr_v(1-r_v)(\widetilde C\eta)_v -\sum_{\{v,w\}}m_{vw}\eta_v\eta_w(r_v-r_w)^2.\end{equation}
Both terms are nonpositive. Equality forces all $r_v$ to equal one
number $r$, since the support is connected. When $\eta$ is primitive,
integrality forces $r$ to be an integer, impossible for a proper
nonzero subvector. When $q_{\Qbar}(\eta)<0$, at least one simple pairing
is strictly negative, so equality is again impossible for $0<r<1$.
Thus $(\gamma,\eta-\gamma)<0$.

Consider the direct-sum morphism
\[ G_\eta\times\Rep(\Qbar,\gamma)\times\Rep(\Qbar,\eta-\gamma) \longrightarrow\Rep(\Qbar,\eta),\qquad (g,X,Y)\longmapsto g\cdot(X\oplus Y). \]
Its fibres contain the free orbits of
$G_\gamma\times G_{\eta-\gamma}$ obtained by changing bases separately
in the two summands and compensating in $g$. Subtracting this group
dimension from the dimension of the source bounds the image closure by
\[ \dim\Rep(\Qbar,\eta)+(\gamma,\eta-\gamma) <\dim\Rep(\Qbar,\eta). \]
There are finitely many proper subvectors, so the complement of all
these closures is a nonempty open subset. This is the direct-sum
argument of \cite[Lemmas~8.3.2--8.3.3]{DWbook};
\eqref{eq:fundamental-identity} gives the strict inequality needed here.
\end{proof}

\subsection{Absolute indecomposability in arbitrary root dimension}

\begin{theorem}[Absolutely indecomposable representations]
\label{thm:absolute-root}
For every nonzero rank vector $\beta$, an absolutely
indecomposable $\Lambda$-module of rank $\beta$ exists if and
only if $\iota\beta$ is a positive root of the unfolded quiver.
If $\iota\beta$ is real, there is exactly one isomorphism class.
If it is imaginary, there is a family $W$ of absolutely indecomposable
representations with
\[ \nu(W)=1-q_d(\beta)>0. \]
These assertions hold over every finite extension of $K$, with
its corresponding partially split species.
\end{theorem}
\begin{proof}
Necessity is ordinary Kac theory after algebraic closure. For a
real root, uniqueness of the geometric indecomposable makes its
isomorphism class Galois-invariant, so orbit descent gives the
unique absolute model. The isotropic case is
Proposition~\ref{prop:isotropic-absolute}.

Suppose $q_d(\beta)<0$. By folding, $\beta$ is a positive
imaginary root. Choose a folded Weyl element carrying it into
the fundamental region. Products of reflections at the pairwise
nonadjacent sheets realize the folded Weyl action. Hence the
resulting unfolded vector $\eta$ is still a positive root and
has connected support and $q_{\Qbar}(\eta)<0$. For every
proper integral subvector $\gamma$, Lemma~\ref{lem:fundamental-open} gives
$(\gamma,\eta-\gamma)<0$.
Removing the proper direct-sum closures therefore leaves a
nonempty open set of geometrically indecomposable representations,
stable under the combined action of Galois on coefficients and
the unfolding. Its form is open in a $K$-affine
space and has a $K$-point. On the open subset where the endomorphism dimension takes its minimum
$h$, we have $\nu=h-q_{\Qbar}(\eta)\ge1-q_{\Qbar}(\eta)$.
Lemma~\ref{lem:orbit-parameter-bound} gives the reverse inequality.
Thus equality holds.

Write the inverse Weyl group element as a product of simple
reflections and apply them in sequence. Before each reflection, use
Lemma~\ref{lem:negative-transport} to choose an orientation in which
the vertex is a sink. Then apply the frame construction of
Lemma~\ref{lem:family-reflection}.
On geometric points this is a product of ordinary reflection
functors, so it preserves absolute indecomposability. No positive
imaginary root is a vertex simple, and therefore no intermediate
module is annihilated. Reflection preserves $\nu$ and the
quadratic form. Finally restore the prescribed orientation by
Lemma~\ref{lem:negative-transport}. This gives the desired family.

After a finite extension of $K$, the split vertices are grouped
into the orbits of a power of $a$. Its vertex factors are still
Laurent-series fields of the same form, so all preceding arguments
apply. A positive-dimensional rational quotient and
Lemma~\ref{lem:henseldense} give infinitely many isomorphism
classes in each imaginary case.
\end{proof}

\begin{corollary}\label{cor:full-generic-moduli}
Let $\beta$ be primitive and let $\theta\in \mc H_\beta$ be
geometrically generic. The moment fibre has a $K$-point if and
only if $\iota\beta$ is a positive root. When nonempty, its
geometric points are simple deformed-preprojective representations,
and its quotient by $G_\beta/\mb G_m$ is smooth of dimension
$2(1-q_d(\beta))$ and has a $K$-point.
\end{corollary}
\begin{proof}
The absolute model supplied by Theorem~\ref{thm:absolute-root}
lifts by the trace criterion. Conversely, genericity forces the
hereditary restriction of a geometric lift to be indecomposable.
It also excludes proper subrepresentations of the doubled
representation. Since every geometric representation in the moment fibre is simple,
its orbit is closed with scalar stabilizer. The affine quotient by
$G_\beta/\mb G_m$ is therefore geometric. The moment-map differential
is surjective onto the annihilator of the scalars, so the fibre and
quotient are smooth; the latter has dimension $2(1-q_d(\beta))$,
as in \cite[Section~2.1]{CBVdB}.
\end{proof}

The parameters in Corollary~\ref{cor:full-generic-moduli} may depend
on $t$. Corollary~\ref{cor:smooth-trace} will concern constant invariant
parameters, for which the cohomology carries the finite diagram action.

\begin{example}[A divisible isotropic dimension]\label{ex:divisible-isotropic}
For $m\ge1$ and $c\in K$, let $J_m(c)$ be a single Jordan matrix
of size $m$ with eigenvalue $c$. Consider the representation $M_c$
of the quiver with two arrows $1\to2$ and one arrow $2\to3$:
\[
\quiverdiagram{\xymatrix@C=3.6em{K^m\ar@/^0.65pc/[r]^{I_m}\ar@/_0.65pc/[r]_{J_m(c)}&K^m\ar[r]^{I_m}&K^m}},\qquad c\in K,\]
The identity maps identify the three endomorphism components, so
an endomorphism is a matrix commuting with $J_m(c)$. Hence
$\End(M_c)\iso K[\varepsilon]/(\varepsilon^m)$.
This algebra remains local after extending the field, and $M_c$ is
absolutely indecomposable. Different values of $c$ give nonisomorphic
representations. Their dimension is $m\gamma$, where
$\gamma=(1,1,1)$ and $q_Q(\gamma)=3-2-1=0$. Thus isotropic
representations at divisible dimensions can occur on an indefinite
support, not only on an affine quiver.
\end{example}

\section{Arithmetic classification by rank and residue degree}
\label{sec:classification}
The preceding section determined the ranks admitting absolutely
indecomposable representations and constructed maximal-dimensional
families at imaginary roots. We now prescribe the Galois-orbit sizes
of their geometric indecomposable summands. At a real root, the isomorphism class has the
same Galois period as its dimension vector. For an imaginary root,
a family can realize any positive multiple of the dimension period.

\subsection{Orbit norms and prescribed periods}

For arbitrary $\beta\ne0$ and $e\geq1$ put
\begin{equation}\label{eq:allrank-root-orbits} \ms R_{\beta,e} =\left\{\alpha\in\Delta^+(\Qbar): a^e\alpha=\alpha,\quad \mc N_e^a\alpha=\iota\beta \right\}/\langle a\rangle.\end{equation}
Write $s(\alpha)$ for the least $a$-period of the dimension
vector, and let $\ms R^{\mr{im}}_{\beta,e}$ be the
imaginary part. Among the real root orbits, we use those of exact period $e$:
\[ \ms R^{\mr{re,exact}}_{\beta,e} =\{[\alpha]\in\ms R_{\beta,e}:\alpha\text{ real},\ s(\alpha)=e\}. \]

\begin{lemma}[Prescribed Galois periods]
\label{lem:allrank-exact-period}
Let $\alpha$ be a positive imaginary root of the unfolding, and let
$s=s(\alpha)$. For every multiple $e$ of $s$, there are $|K|$ Galois
orbits of indecomposable $\Lambda_{\Kbar}$-modules that contain a module
of dimension $\alpha$ and have size $e$. The corresponding
$\Lambda$-modules have rank $f(\mc N_e^a\alpha)$ and residue
endomorphism field $K_e$.
\end{lemma}
\begin{proof}
Put $F=K_s$, $E=K_e$, and $m=e/s$. Write
$\OO_F=k[[t^{1/s}]]$ for the valuation ring of $F$.
Theorem~\ref{thm:absolute-root}
for the partial species over $F$ supplies a family of absolutely indecomposable representations of
dimension $\alpha$ with positive $\nu$. Rosenlicht's theorem gives
a nonconstant invariant rational function $\varphi$ defined over $F$.
After shrinking, $\varphi$ is regular and has nonzero differential.
Density of rational points supplies an $F$-point $z_0$ where it is
smooth. Write $z=t^{1/s}$ and $v=t^{1/e}$, so $v^m=z$.
Hensel's lemma over $E$ realizes the values
\[ \varphi(z_0)+z^nvc, \qquad c\in\OO_F^\times, \]
for one sufficiently large $n$, independent of the unit $c$.
Indeed all these increments have the same sufficiently large valuation.
Each value generates $E$ over $F$, since $c z^n$ is a nonzero element
of $F$ and $F(v)=E$. The invariant $\varphi$ is defined over $F$; it is
constant on geometric isomorphism classes and satisfies
$\varphi(\sigma M)=\sigma \varphi(M)$ for $\sigma\in\Gal(E/F)$.
Thus the underlying absolute $E$-module has $m$ distinct
conjugates within the stabilizer of its dimension vector. A
conjugation outside that stabilizer changes the dimension vector.
The full orbit therefore has exactly $e$ distinct members.
Orbit descent, or restriction of scalars from $E$ to $K$, gives
the asserted indecomposable and residue field.

The set of units of $\OO_F$ has cardinality $|K|$.
These values are distinct, and each full conjugacy orbit
identifies at most $e$ of the resulting isomorphism classes.
Thus $|K|$ classes remain. The upper bound $|K|$ follows from
the finite matrix description of representations.
\end{proof}

\begin{theorem}[Arithmetic classification]
\label{thm:allrank-strata}
For every nonzero rank $\beta$ and every $e\geq1$,
\[
\#\Ind_{\beta,e}(\Lambda)=
 \begin{cases}
   |K|,&\ms R^{\mr{im}}_{\beta,e}\ne\varnothing,\\[2pt]
   |\ms R^{\mr{re,exact}}_{\beta,e}|,
      &\ms R^{\mr{im}}_{\beta,e}=\varnothing.
 \end{cases}
\]
The root-orbit conditions give a finite test for the cardinality
of $\Ind_{\beta,e}(\Lambda)$. It depends on the chosen modulation but
not on the acyclic orientation. The set is empty unless
$e\mid g_d(\beta)$.
\end{theorem}
\begin{proof}
An indecomposable with residue field $K_e$ has one multiplicity-free
geometric orbit of $e$ indecomposables. Its dimension is a root
$\alpha$ in \eqref{eq:allrank-root-orbits}. If $\alpha$ is real,
uniqueness of the geometric representation forces its isomorphism
period to equal $s(\alpha)$, hence $s(\alpha)=e$. Conversely each
such real orbit has exactly one descent. An imaginary candidate
supplies $|K|$ descents by Lemma~\ref{lem:allrank-exact-period}.
This proves the alternatives. The divisibility follows also from
Proposition~\ref{prop:divisibility}. All root and norm
conditions depend on the unoriented unfolding and its action,
not the orientation.
\end{proof}

The bound $0\le\alpha\le\iota\beta$ leaves only finitely many
solutions of the norm equations. Standard height reduction decides
which solutions are roots; among real roots, only the orbits of exact
period $e$ are retained.

\subsection{Kac-type root classification and residue-field towers}
We now deduce the root classification for the species and exhibit
residue degrees occurring at every positive multiple of an imaginary
root. The following lemma constructs the required unfolded roots.

\begin{lemma}\label{lem:imaginary-preimage}
For every positive imaginary root $\beta$ of $C$ there is a positive
imaginary root $\alpha$ of the unfolding, of least period $s\mid N$,
such that $\mc N_s^a\alpha=\iota\beta$.
\end{lemma}
\begin{proof}
We choose $w$ in the folded Weyl group such that $\gamma=w\beta$ lies
in the fundamental region. Restrict the unfolding to the inverse
image of $\supp\gamma$. Its connected components are permuted
transitively by $a$: the translates of one component contain whole
vertex orbits and all their neighbours. Let $c_\gamma$ be the number of these connected components.
The restriction $\eta$ of $\iota\gamma$ to one component, extended
by zero, has connected support and all simple pairings nonpositive.
It is therefore an imaginary root. It has least period $c_\gamma$, and
$\mc N_{c_\gamma}^a\eta=\iota\gamma$.

By \cite[Lemma~3]{Hubery}, a product of the commuting simple
reflections over a vertex orbit realizes a folded simple reflection
and commutes with $a$.
Lift $w^{-1}$ in this way and apply it to $\eta$, obtaining $\alpha$.
Positive imaginary roots remain positive under the Weyl group.
The lift is invertible and commutes with $a$, so it preserves the
least period. Its orbit norm is $\iota\beta$, as required.
\end{proof}

\begin{corollary}\label{cor:full-kac}
The indecomposable ranks are exactly $\Delta^+(C)$. For a real root,
the unique class and its residue field are described in
Proposition~\ref{prop:root}; an imaginary root has $|K|$ classes. For every positive imaginary root $\beta$, some
$s\mid\gcd(N,g_d(\beta))$ satisfies
\[ \#\Ind_{m\beta,sm}(\Lambda)=|K|\qquad(m\ge1). \]
One may take $s=1$ when $g_d(\beta)=1$. For the primitive null root
$\delta$ of a connected affine $C$, one may take $s=\unfoldcomp$, the number
of unfolded connected components.
\end{corollary}
\begin{proof}
Use the imaginary preimage of Lemma~\ref{lem:imaginary-preimage}.
Its norm at period $sm$ is $\iota(m\beta)$, so
Theorem~\ref{thm:allrank-strata} applies; the divisibility follows
from Proposition~\ref{prop:divisibility}. In affine type $C\delta=0$,
so restriction to a component already gives period $\unfoldcomp$, without
applying any reflections. Proposition~\ref{prop:root} gives the set
of ranks and the real-root assertion.
\end{proof}

\section{Applications for affine types}\label{sec:affine-strata}
Let $C$ be connected affine. Its \emph{primitive null root}
$\delta$ is the positive integral generator of $\ker C$.
We use the classical tube classification of ordinary affine quivers
\cite[Section~7.3]{DWbook}, \cite{DlabRingel}, to make the arithmetic
criterion explicit. Both tube ranks and the action on their
quasi-simples enter the answer.

The $\unfoldcomp$ connected components of $\Qbar$ are permuted transitively by $a$.
Choose one, $Q^\circ$. Its return automorphism is $a^\unfoldcomp|_{Q^\circ}$,
and its field is $F=K_\unfoldcomp$. The restriction $\delta_0=(\iota\delta)|_{Q^\circ}$ is a primitive
positive null vector, so $Q^\circ$ is affine. Its \emph{homogeneous tubes} are parametrized by
$\mb P^1_k\setminus\Sigma$, where the finite set $\Sigma$
indexes the \emph{exceptional tubes}. Descent gives an open curve $U_F$ whose
completion is $\mb P^1_F$, since $\Br(F)=0$.

For $Q^\circ$, the \emph{Euler defect} of a representation $X$ is
$\partial(X)=\langle\delta_0,\dimv X\rangle_{Q^\circ}$, using the
ordinary Euler form. The indecomposables of defect zero are called
\emph{regular}. The regular category decomposes into extension-closed
abelian subcategories called \emph{tubes}. A tube of rank $r$ has $r$
simple objects, called its \emph{quasi-simples}. Each indecomposable in
the tube has a unique simple regular submodule, its
\emph{quasi-socle}, and a filtration by quasi-simples. The length of
this filtration is its \emph{quasi-length}. Such an indecomposable
is \emph{uniserial}: its regular submodules form a chain. A homogeneous tube has
rank one; the exceptional tubes have rank greater than one.

For an orbit $\OO$ of exceptional tubes under $a^\unfoldcomp$, let $h_{\OO}$
be its size and $r_{\OO}>1$ the tube rank. The return map
$a^{\unfoldcomp h_{\OO}}$ shifts the cyclically indexed quasi-simples by
$s_{\OO}$. Define
\[ g_{\OO}=\gcd(r_{\OO},s_{\OO}),\qquad e_{\OO}=\unfoldcomp h_{\OO}r_{\OO}/g_{\OO}, \qquad\gcd(r,0)=r. \]
For a closed point $x$ on a curve over $F$, its degree means
$[\kappa(x):F]$, where $\kappa(x)$ is its residue field. The omitted
closed points of $U_F$ have degrees $h_{\OO}$.
A cyclic shift of the starting index leaves $s_{\OO}$ unchanged,
whereas reversing the cyclic indexing replaces it by $-s_{\OO}$.
In either case $g_{\OO}$ and $e_{\OO}$ are unchanged.

\begin{proposition}\label{prop:affine-strata}
The indecomposables of rank $m\delta$ are the disjoint union of:
\begin{enumerate}
\item one class from a homogeneous tube for each closed point $x\in U_F$ of degree
$d\mid m$, with geometric quasi-length $m/d$ and residue field $K_{\unfoldcomp d}$;
\item $g_{\OO}$ exceptional classes for each orbit $\OO$
with $h_{\OO}\mid m$, with geometric quasi-length
$g_{\OO}m/h_{\OO}$ and residue field $K_{e_{\OO}}$.
\end{enumerate}
Consequently,
\[ \Ind_{m\delta,e}(\Lambda)\text{ is infinite} \quad\Longleftrightarrow\quad e=\unfoldcomp d\text{ for some }d\mid m. \]
For every other $e$ its cardinality is
\[ \#\Ind_{m\delta,e}(\Lambda) =\sum_{\substack{\OO:\ h_{\OO}\mid m\\e_{\OO}=e}} g_{\OO}. \]
\end{proposition}
\begin{proof}
A geometric summand of a module of rank $m\delta$ is regular:
all conjugates have the same sign of Euler defect, and their total
defect is zero. A homogeneous summand is determined by its parameter
and quasi-length $n$. A parameter of degree $d$ on $U_F$ has full
Galois orbit $\unfoldcomp d$ and norm $nd\,\iota\delta$. Orbit descent gives
exactly (i).

For an exceptional-tube orbit, abbreviate its data by $h,r,s,g$. The quasi-socle
and length determine an exceptional indecomposable, including when
$r$ divides its length. At each length $n$, addition by $s$ on the
$r$ quasi-socles gives $g$ module orbits, each of full size $\unfoldcomp hr/g$.
Let $\zeta_j$ be their quasi-simple dimensions. They satisfy
\[ \sum_{j=0}^{r-1}\zeta_j=\delta_0,\qquad \langle\zeta_i,\zeta_j\rangle =\delta_{ij}-\delta_{j,i+1}, \]
with cyclic indices; distinct tubes are Euler-orthogonal.
The sum of a length-$n$ cyclic interval over the return orbit has
constant coefficients on the $\zeta_j$ exactly when $g\mid n$.
Indeed these coefficients count positions in each residue class modulo
$g$. Pairing with $\zeta_i$ detects their consecutive differences,
so nonconstant contributions cannot cancel between tubes. A constant
sum contributes $(n/g)\delta_0$ from each of the $h$ tubes. Thus its
norm is $m\iota\delta$ exactly when $m=hn/g$.
Coefficient conjugation fixes these constant-field modules; their
Galois orbits are their $a$-orbits. This proves (ii).

Finally, $U_F$ has infinitely many closed points of every positive
degree: after writing $F=k((z))$, the polynomials $X^d-\eta z$,
$\eta\in k^\times$, give distinct degree-$d$ points for $d>1$.
Deleting finitely many points changes none of these infinitude assertions.
\end{proof}

\subsection{Finite-field counts}
For the species $\Lambda_{q,a}$ defined in the Introduction, the same tube
orbits determine the counts at every imaginary rank. Put
\[ N_d(T)=\frac1d\sum_{j\mid d}\mu_{\mr{ar}}(d/j)(T^j+1), \qquad b_d=\#\{\OO:h_{\OO}=d\}. \]
At a prime power $T$, $N_d(T)$ counts degree-$d$ closed points of
$\mb P^1_{\mb F_T}$; in particular $N_1(0)=1$ and $N_d(0)=0$
for $d>1$.

\begin{corollary}\label{cor:affine-counts}
For connected affine $C$, every $m,e\geq1$, and every finite field,
\begin{equation}\label{eq:affine-counts}
 I_{m\delta,e}(q)
 =\begin{cases}
   N_d(q^\unfoldcomp)-b_d,&e=\unfoldcomp d\text{ with }d\mid m,\\
   0,&\text{otherwise}
  \end{cases}
  +\sum_{\substack{\OO:\ h_{\OO}\mid m\\
                                e_{\OO}=e}}g_{\OO}.
\end{equation}
In particular the total count is the integer-valued rational
polynomial
\begin{equation}\label{eq:affine-total-count} I_{m\delta,a}(q) =\sum_{d\mid m}N_d(q^\unfoldcomp) +\sum_{\substack{\OO:h_{\OO}\mid m}} (g_{\OO}-1).\end{equation}
Its constant coefficient is
\begin{equation}\label{eq:affine-constant-count} I_{m\delta,a}(0)=1+ \sum_{\substack{\OO:h_{\OO}\mid m}} (g_{\OO}-1)\ \geq1.\end{equation}
\end{corollary}
\begin{proof}
Over $\mb F_{q^\unfoldcomp}$ the homogeneous completion is a form of
$\mb P^1$, hence is split. Ordinary Frobenius fixes the exceptional
quasi-simples and their uniserial extensions; the quasi-simple real roots
have unique representations with prime-field models. Thus the omitted
closed points have degrees $h_{\OO}$, and the exceptional orbit
calculation in Proposition~\ref{prop:affine-strata} is unchanged, in every
characteristic. There are $N_d(q^\unfoldcomp)-b_d$ degree-$d$ homogeneous points.
Combining the two lists proves \eqref{eq:affine-counts}. Summing over
$e$ gives \eqref{eq:affine-total-count}; $N_1(0)=1$ and $N_d(0)=0$
for $d>1$ give \eqref{eq:affine-constant-count}.
\end{proof}

\subsection{Two arithmetic distinctions}
\begin{example}\label{ex:two-modulations}
Let $L=K_2$. The bimodules $L\oplus L$ and $L\oplus L^{\mr{tw}}$,
where one field action is twisted by the nontrivial element of
$\Gal(L/K)$, both give
$C=\bigl(\begin{smallmatrix}2&-2\\-2&2\end{smallmatrix}\bigr)$,
$D=\op{diag}(2,2)$ and $\delta=(1,1)$. Choose the
orientation from vertex $1$ to vertex $2$. The two species and their
unfoldings are
\[
\begin{array}{c@{\qquad\qquad}c}
\quiverdiagram{\xymatrix@C=3em{L\ar@/^0.6pc/[r]^{L}\ar@/_0.6pc/[r]_{L}&L}}
&\quiverdiagram{\xymatrix@C=3em{L\ar@/^0.6pc/[r]^{L}\ar@/_0.6pc/[r]_{L^{\mr{tw}}}&L}}\\[1.8em]
\quiverdiagram{\xymatrix@C=3em@R=2em{
u_0\ar@/^0.55pc/[r]^{b_0}\ar@/_0.55pc/[r]_{c_0}&v_0\\
u_1\ar@/^0.55pc/[r]^{b_1}\ar@/_0.55pc/[r]_{c_1}&v_1}}
&\quiverdiagram{\xymatrix@C=2.3em@R=1.7em{
 &u_0\ar[dl]_{b_0}\ar[dr]^{c_0}&\\
 v_0&&v_1\\
 &u_1\ar[ul]^{c_1}\ar[ur]_{b_1}&}}
\end{array}
\]

In both drawings $a$ exchanges $u_0,u_1$ and $v_0,v_1$; it also
exchanges the two arrows of each labelled type. The first unfolding
is two Kronecker components without exceptional
tubes. At $m\delta$ precisely the residue degrees $2d$, with $d\mid m$, occur, and each of these
sets is infinite.
In particular, at $2\delta$ residue degree four occurs although
$N=2$. A field extension that splits the algebra need not split
every indecomposable representation.
The second unfolding is the connected square $\mathsf K_{2,2}$.
Its two rank-two exceptional tubes are fixed, with quasi-simples
exchanged. For each residue degree dividing $m$, the corresponding set is infinite; for odd $m$ there
are also two degree-two classes, and all other residue degrees give empty sets.

Over a finite field, at $\delta$ the untwisted species has
$I_{\delta,2}(T)=T^2+1$ and $I_{\delta,1}(T)=0$; the twisted species
has $I_{\delta,1}(T)=T-1$ and $I_{\delta,2}(T)=2$.
Thus neither the distribution of residue degrees nor the
refined counting polynomials are determined by $C$ and $D$ alone.
\end{example}

\begin{example}\label{ex:14-affine}
For $C=\bigl(\begin{smallmatrix}2&-4\\-1&2\end{smallmatrix}\bigr)$,
$D=\op{diag}(1,4)$ and $\delta=(2,1)$, the unfolding is the
four-leaf $\widetilde D_4$ star \eqref{eq:1,4}.
The three rank-two tubes correspond to partitions of the leaves into
pairs. One is fixed and two are exchanged; each return swaps the
quasi-simples. Thus
\[ (h,r,s,g,e)=(1,2,1,1,2),\qquad(2,2,1,1,4). \]
At rank $m\delta$, the sets whose residue degrees divide $m$ are
infinite. Outside them there
is one degree-two class for odd $m$, one degree-four class for
$m\equiv2\pmod4$, and no other classes; see also the matrix
constructions in \cite{GR,GRR}.
The total finite-field count is $\sum_{d\mid m}N_d(q)$, as for the
split Kronecker quiver, but the refined counts differ:
\[
\begin{array}{c|ccc}
 &e=1&e=2&e=4\\\hline
 m=1&q&1&0\\
 m=2&q&\tfrac12(q^2-q)&1
\end{array}
\]
For the split Kronecker quiver at $2\delta$, the counts in residue
degrees $1,2,4$ are instead $q+1$, $\tfrac12(q^2-q)$, and $0$. Agreement of total counts does not
identify the sets with prescribed residue degree.
\end{example}

\section{Finite-field counts and nilpotent endomorphisms}
\label{sec:count}

We now fix the unfolding $(\Qbar,a)$ of
Proposition~\ref{prop:unfold}, viewed as combinatorial data over
$\mb Z$. For a finite field $\mb F_q$, let
\[ \Lambda_{q,a}=(\overline{\mb F}_q\Qbar)^{F_a},\qquad F_a=a\circ\mr{Frob}_q, \]
where $F_a$ acts semilinearly on coefficients and permutes vertices and
arrows. This is a finite-dimensional species over $\mb F_q$,
with vertex fields $\mb F_{q^{d_i}}$. The simple factors of its
arrow bimodules are specified by $a$. A representation of rank
$\beta$ has dimension vector $\iota\beta$ over
$\overline{\mb F}_q$. We also write $V_\beta$ and $G_\beta$
for the representation space and change-of-basis group over this
finite field.

Extending $\mb F_q$ to $\mb F_{q^r}$ replaces the descent
operator by $a^r\circ\mr{Frob}_{q^r}$; it does not keep $a$ fixed.
When $a$ is held fixed below, the species is reconstructed for each $q$.
These facts are the fixed-algebra construction and scalar-extension
formula of \cite[Proposition~6.3 and Corollary~6.4]{DengDu}.

The root correspondence is known from \cite{DengXiao}; polynomiality
and orientation independence of the total indecomposable count are
\cite[Corollary~9.2]{DengDu}. We now separate this count by residue
endomorphism fields and compare it with the action on cohomology.

Let $I_{\beta,e}(q)$ count the indecomposables of rank $\beta$
whose residue endomorphism field is $\mb F_{q^e}$. By
\cite[Theorem~5.1 and Corollary~5.2]{DengDu}, its geometric summands
form one Frobenius orbit of length $e$, each occurring once; compare
\cite[Lemma~8.4.6]{DWbook} and \cite[Lemma~19]{Hubery}.
Consequently $e\mid g_d(\beta)$, by the proof of
Proposition~\ref{prop:divisibility}. We omit $a$ from this notation
when the unfolding is fixed.

\subsection{Exact residue degrees under partial splitting}

Keep $F_a=a\circ\mr{Frob}_q$ fixed. For a nonnegative unfolded
vector $\alpha$ and $r\geq1$, let $\mathsf A_{\alpha,r}(q)$ be the number
of geometric indecomposable isomorphism classes of dimension
$\alpha$ fixed by $F_a^r$. We set this number to zero when
$a^r\alpha\ne\alpha$. When $a^r\alpha=\alpha$, it is the number
of absolutely indecomposable representations, with that split
dimension vector, of the species descended by
$a^r\circ\mr{Frob}_{q^r}$ over $\mb F_{q^r}$.
These descriptions agree by \cite[Corollary~5.4]{DengDu}.
Existence of a Frobenius structure is the descent assertion; uniqueness
of the rational isomorphism class follows from
\cite[Lemma~3.1]{DengDu}, or from Lang's theorem for the connected
unit group of the geometric endomorphism algebra.
We use the orbit norm $\mc N_e^a$ and the arithmetic
M\"obius function $\mu_{\mr{ar}}$ defined in the Introduction.

\begin{lemma}\label{lem:mobius}
For every finite field and every $e\geq1$, we have
\begin{equation}\label{eq:mobius-residue} I_{\beta,e}(q)=\frac1e \sum_{\substack{\alpha\in\mb N^{\Qbar_0}\setminus\{0\}\\ \mc N_e^a\alpha=\iota\beta}} \ \sum_{r\mid e}\mu_{\mr{ar}}(e/r)\mathsf A_{\alpha,r}(q).\end{equation}
The outer sum is finite. It is empty unless $e\mid g_d(\beta)$.
\end{lemma}
\begin{proof}
For fixed $\alpha$, the number of geometric indecomposable classes
of dimension $\alpha$ having exact $F_a$-period $e$ is
$\sum_{r\mid e}\mu_{\mr{ar}}(e/r)\mathsf A_{\alpha,r}(q)$.
Indeed, for each class fixed by $F_a^e$, the sum of
$\mu_{\mr{ar}}(e/r)$ over the divisors $r$ for which $F_a^r$
fixes it is one if its least period is $e$, and zero otherwise.
This argument does not require $F_a$ to preserve $\alpha$; terms
with $a^r\alpha\ne\alpha$ are simply zero.

An orbit of exact period $e$ descends to an indecomposable whose
residue endomorphism field has degree $e$. Its total geometric
dimension is $\mc N_e^a\alpha$. Conversely, every
indecomposable of residue degree $e$ is obtained this way. Summing
over all possible $\alpha$ counts every such orbit exactly $e$
times, which proves \eqref{eq:mobius-residue}. The norm equation
implies $\alpha\leq\iota\beta$, so the sum is finite.
It also implies $a^e\alpha=\alpha$, since
$(a-1)\mc N_e^a\alpha=a^e\alpha-\alpha$.
Summing the norm equation over the sheets above $i$ gives
$d_i\beta_i=e\sum_{v\text{ above }i}\alpha_v$,
which proves the asserted divisibility.
\end{proof}

For $a=\mr{id}$, the formula reduces to the ordinary identity
\[ I_{\beta,e}(q)=\frac1e\sum_{r\mid e} \mu_{\mr{ar}}(e/r)A_{\beta/e}(q^r) \quad(e\mid\gcd_i\beta_i), \]
and is zero otherwise; compare \cite[Lemma~8.4.7]{DWbook}.
For nontrivial $a$, the counts on the right belong to the partial
splittings defined by $a^r$. Theorem~\ref{thm:allrank-absolute} will
identify each of these absolute counts with a cohomological trace.

\begin{remark}\label{rem:nonintegral-count}
For the split Kronecker quiver with $\delta=(1,1)$, the ordinary
Kac polynomials satisfy $A_{Q,\delta}(T)=A_{Q,2\delta}(T)=T+1$.
Lemma~\ref{lem:mobius} gives
\[ I_{2\delta,2}(T)=\frac{(T^2+1)-(T+1)}2=\frac{T^2-T}{2},\qquad
   I_{2\delta,\mr{id}}(T)=\frac{T^2+T+2}{2}. \]
The half-integral coefficients arise from counting Galois orbits
of length two, even for split quivers.
\end{remark}

\subsection{Nilpotent endomorphisms and a Jordan formula}

For a variety $Z$ with an algebraic group action, $[Z/G]$ denotes
the quotient stack. For a stack $\mathfrak X$ over $\mb F_q$,
its \emph{stack cardinality} is
$|\mathfrak X(\mb F_q)|=\sum_{[x]}1/|\Aut(x)|$.
Let $\mathbf x=(x_i)_{i\in I}$ be formal variables and put
$\mathbf x^\beta=\prod_i x_i^{\beta_i}$. For $\Lambda_{q,a}$,
consider the \emph{stack of representations with a nilpotent
endomorphism} and its counting series:
\[ \mathfrak N_{\beta,a} =\bigl[\{(M,f):M\in V_\beta,\ f\in\End(M),\ f\text{ nilpotent}\}/G_\beta\bigr], \qquad \mc J_a(\mathbf x;q) =\sum_{\beta\geq0}|\mathfrak N_{\beta,a}(\mb F_q)|\,\mathbf x^\beta. \]
Here $\End(M)$ is taken over $\Lambda_{q,a}$; the same stack
notation without $a$ denotes the split complex stack at an unfolded
dimension vector.
We write $\lambda\vdash m$ for a partition of $m$,
$\lambda'_r=\#\{j:\lambda_j\geq r\}$ for the parts of its
conjugate partition, and $m_s(\lambda)=\#\{j:\lambda_j=s\}$.
Empty partitions are allowed. Put
\begin{equation}\label{eq:jordan-coefficient} J_{\beta,a}(T)= \sum_{\lambda^i\vdash\beta_i} \frac{T^{-\sum_{r\geq1}q_d(\mathbf c_r)}} {\displaystyle\prod_i\prod_{s\geq1} \prod_{j=1}^{m_s(\lambda^i)}(1-T^{-d_i j})}, \qquad \mathbf c_r=((\lambda^i)'_r)_i.\end{equation}
This is a finite sum in $\mb Q(T)$, and $J_{0,a}=1$.

\begin{proposition}\label{prop:nil-jordan}
For every finite field, we have
\begin{align}
 |\mathfrak N_{\beta,a}(\mb F_q)|&=J_{\beta,a}(q),\label{eq:jordan-count}\\
 \mc J_a(\mathbf x;q)
 &=\prod_{\gamma\ne0}\prod_{e\mid g_d(\gamma)}
       \left(\prod_{j\geq1}(1-q^{-ej}\mathbf x^\gamma)^{-1}\right)^{I_{\gamma,e}(q)}.
       \label{eq:nil-product}
\end{align}
For the fixed automorphism $a$, the expression
\eqref{eq:jordan-coefficient} depends only on $C$ and $D$.
\end{proposition}
\begin{proof}
We fix the Jordan type $\lambda^i$ of $f_i$ over $\mb F_{q^{d_i}}$.
The space of arrow maps commuting with the $f_i$ has dimension
\[ \sum_{i,j}(\dim_K A_{ij})\sum_{r\geq1}(\lambda^i)'_r(\lambda^j)'_r \]
over $\mb F_q$. Here $\dim_K A_{ij}$ also equals the base-field
dimension of the corresponding finite-field arrow bimodule. The formula
holds for every simple bimodule factor: after scalar extension it is
the usual identity $\dim\Hom(k[z]/z^u,k[z]/z^v)=\min(u,v)$.
The centralizer of $f_i$ has dimension
$d_i\sum_r((\lambda^i)'_r)^2$ and reductive quotient
$\prod_s\Res_{\mb F_{q^{d_i}}/\mb F_q}
\mr{GL}_{m_s(\lambda^i)}$. Its number of units is therefore
\[ q^{d_i\sum_r((\lambda^i)'_r)^2} \prod_s\prod_{j=1}^{m_s(\lambda^i)}(1-q^{-d_i j}). \]
The exponent of $q$ after division is
$-\sum_r\langle\mathbf c_r,\mathbf c_r\rangle_A
=-\sum_r q_d(\mathbf c_r)$, which explains the exponent in
\eqref{eq:jordan-coefficient}.
Dividing the number of compatible arrow maps by these centralizer
orders and summing proves \eqref{eq:jordan-count}. This is Hua's Jordan-form calculation \cite{Hua}, with the
vertex-field degrees included.

For the second assertion write $M=\bigoplus_X X^{n_X}$ and let
$\End(X)/\rad\End(X)=\mb F_{q^{e_X}}$. An element of $\End(M)$ is
nilpotent exactly when its image in the semisimple quotient is
nilpotent. We write $\op{Nil}\End(M)$ for the set of
nilpotent endomorphisms. For a prime power $Q$, in
$\op{Mat}_n(\mb F_Q)$ there are
$Q^{n^2-n}$ nilpotent matrices. Each nilpotent element of the semisimple quotient has
$q^{\dim\rad\End(M)}$ lifts, and the same factor occurs in
$|\Aut(M)|$. It cancels, giving
\[ \frac{|\op{Nil}\End(M)|}{|\Aut(M)|} =\prod_X\frac{q^{-e_X n_X}} {\prod_{j=1}^{n_X}(1-q^{-e_Xj})}. \]
We now use the $q$-binomial identity
\[ \sum_{n\geq0}\frac{Q^{-n}z^n}{\prod_{j=1}^n(1-Q^{-j})} =\prod_{j\geq1}(1-Q^{-j}z)^{-1} \]
Applying it independently to the multiplicities $n_X$ proves
\eqref{eq:nil-product}.
\end{proof}

Although \eqref{eq:jordan-coefficient} depends only on $C,D$, the
individual residue-degree counts use this formula for all powers of
$a$. Their partial splittings can differ for two modulations with the
same original Cartan matrix and symmetrizer; see
Example~\ref{ex:two-modulations}.

\subsection{The preprojective counting series}
Let $\mathfrak P_\beta$ be the stack of representations of the
preprojective algebra, and let $\mathfrak E_\beta$ be the stack of
ordinary representations with an arbitrary endomorphism:
\[ \mathfrak P_\beta=[\mu_\beta^{-1}(0)/G_\beta],\qquad \mathfrak E_\beta= [\{(M,f):M\in V_\beta,\ f\in\End(M)\}/G_\beta]. \]
Unlike $\mathfrak N_{\beta,a}$, the second stack allows arbitrary
endomorphisms. Set
\[ \mc Z_a(\mathbf x;q)=\sum_{\beta\ge0} q^{q_d(\beta)}|\mathfrak P_\beta(\mb F_q)|\,\mathbf x^\beta. \]
\begin{corollary}\label{cor:mass-product}
Over every finite field, we have
\[ \mc Z_a(\mathbf x;q)= \prod_{\gamma\ne0}\prod_{e\mid g_d(\gamma)} \left(\prod_{j\ge0}(1-q^{-ej}\mathbf x^\gamma)^{-1}\right)^{I_{\gamma,e}(q)}. \]
\end{corollary}
\begin{proof}
The Euler identity also holds over finite fields. The zero-moment
fibre over $M$ has $q^{h_M-q_d(\beta)}$ points; multiplying by
$q^{q_d(\beta)}$ shows that the coefficient of $\mathbf x^\beta$ is
\[ \sum_{[M]:\rk M=\beta}\frac{q^{h_M}}{|\Aut(M)|}. \]
For $M=\bigoplus_X X^{n_X}$, let $e_X$ be the residue degree of $X$.
The radical cancels and the summand is
$\prod_X\prod_{j=1}^{n_X}(1-q^{-e_Xj})^{-1}$.
The same $q$-binomial identity as in Proposition~\ref{prop:nil-jordan},
now without the factor $q^{-e_Xn_X}$, gives the product starting at
$j=0$. Both products are interpreted coefficientwise for $q>1$.
\end{proof}
The same relation holds in the Grothendieck ring of stacks, in which
we impose additivity for a closed substack and its complement and the
vector-bundle relation. Let $\mb L$ be the class of the affine
line and invert it. Then
\[ [\mathfrak E_\beta]=\mb L^{q_d(\beta)}[\mathfrak P_\beta]. \]
Indeed, on a stratum where $h_M$ is constant, the two bundles have
ranks $h_M$ and $h_M-q_d(\beta)$.
The vector-bundle relation gives the identity,
also equivariantly for the split cyclic action.

\section{BPS cohomology, traces, and root multiplicities}
\label{sec:primitive}\label{sec:allrank-cohomology}
Section~\ref{sec:count} counted representations with a nilpotent
endomorphism by Jordan normal form. We now use equivariant PBW to
express their cohomology in terms of preprojective BPS cohomology.
The stacks $\mathfrak E_\alpha$ and $\mathfrak P_\alpha$, of representations
with unrestricted endomorphisms and preprojective representations,
respectively, enter through dimensional reduction. The degree-zero
BPS action and the twining theorem then give the root multiplicities.

\subsection{Cohomological normalization}
Throughout the cohomological discussion, set $k=\mb C$ and
$K=\mb C((t))$. We work on the unfolded quiver, initially at
$\alpha=\iota\beta$. For an ordinary quiver $Q$, we write
$\Pi_Q=\mb C Q^{\mr{dbl}}/(\rho_v)_{v\in Q_0}$ for its
ordinary preprojective algebra, with the relations of
Section~\ref{sec:moment}.
Write $V_\alpha=\Rep(\Qbar,\alpha)$ and
$G_\alpha=\prod_v\mr{GL}_{\alpha_v}$, and use
$\mathfrak P_\alpha,\mathfrak E_\alpha,\mathfrak N_\alpha$ for the
corresponding split stacks. We write $H_c^i(-,\mb Q)$ for compactly supported rational
cohomology and $H_c$ for its cohomologically graded object. For quotient
stacks we use finite-dimensional approximations. In this section
$\mb L=\mb Q(-1)[-2]$ denotes the \emph{Tate object}, rather
than its class in a Grothendieck ring. The Tate Hodge structure
$\mb Q(-1)$ has weight two, and $[-2]$ places it in cohomological
degree two. A cohomology object is \emph{pure of Tate type} when its
odd pieces vanish and each degree $2r$ is a sum of copies of
$\mb Q(-r)$. We write $H(i)=H\otimes\mb Q(i)$ for a Tate
twist. Tensoring with $\mb L^r$ shifts
cohomological degree and weight by $2r$. At a fixed dimension vector,
the cohomology is finite-dimensional in each degree and vanishes in
sufficiently large degrees; negative degrees may occur. The associated
series are therefore Laurent series in $T^{-1}$. All completions below are taken separately in each dimension vector
and cohomological degree. For an algebraic group $G$, its classifying
stack is $BG=[\op{Spec}\mb C/G]$. When we discuss inertia, we use geometric $\ell$-adic cohomology with
coefficients in $\mb Q_\ell$, for a fixed prime $\ell$, and compare it
with rational cohomology over $\mb C$ after tensoring by
$\mb Q_\ell$. The inertia subgroup of
$\Gal(\overline{\mb C((t))}/\mb C((t)))$ is the whole Galois
group, since the residue field is algebraically closed.

\begin{proposition}[Equivariant dimensional reduction]
\label{prop:allrank-purity}
For every $\beta$, there is an $a$-equivariant isomorphism
\begin{equation}\label{eq:allrank-dimred} H_c(\mathfrak E_\alpha)\iso H_c(\mathfrak P_\alpha)\otimes\mb L^{q_d(\beta)}.\end{equation}
Both sides are pure of Tate type. After the preceding comparison with
$\mb Q_\ell$-cohomology, the action of $a$ agrees with that of a
chosen generator of inertia on the corresponding ramified stacks,
up to inversion of the generator.
\end{proposition}
\begin{proof}
Put $v=\dim V_\alpha$, the dimension of the space of arrow maps.
On $V_\alpha\times V_\alpha^*\times\Lie G_\alpha$ consider
\[ \Phi(M,\xi,f)=\xi(d_Mf)=\mu(M,\xi)(f). \]
This is the trace of the canonical potential on the \emph{tripled quiver}, obtained from $\Qbar$ by adjoining reverse arrows and one
loop at every vertex. Over a single arrow $b:i\to j$, the added
arrows and loops are
\[
\quiverdiagram{\xymatrix@C=4.5em{
i\ar@/^0.65pc/[r]^{b}\ar@(ul,dl)[]_{\omega_i}
&j\ar@/^0.65pc/[l]^{b^*}\ar@(ur,dr)[]^{\omega_j}}}
\]

There is one loop at each vertex, shared by all incident arrows.
Dimensional reduction along $V_\alpha^*$ gives
$H_c(\mathfrak E_\alpha)\otimes\mb L^v$, whereas reduction along
$\Lie G_\alpha$ gives $H_c(\mathfrak P_\alpha)\otimes\mb L^{\dim G_\alpha}$ from the
same vanishing-cycle cohomology
\cite[Section~3.3 and Theorem~3.7]{DavisonStacks}. Both projections
and $\Phi$ are preserved by $a$. Naturality of the reduction maps
therefore gives equivariant isomorphisms. Cancelling $\mb L^v$
proves \eqref{eq:allrank-dimred}, since
$\dim G_\alpha-v=\chi_{\Qbar}(\alpha,\alpha)
=\langle\beta,\beta\rangle_A=q_d(\beta)$.
Davison's Theorem~A gives purity and Tate type for
$H_c(\mathfrak P_\alpha)$, and the same follows for
$H_c(\mathfrak E_\alpha)$. For the assertion about inertia, extend
scalars to $K_N$. Each ramified stack then becomes the base change
of the corresponding complex stack, with descent datum combining
coefficient conjugation and $a$. The coefficient action on its
geometric cohomology is trivial because this model is defined over
$\mb C$; the Tate twists are trivial as Galois representations
since $\mb C$ contains all roots of unity. The remaining action is
therefore the finite action of $a$. Comparison over $\mb C$ gives
the asserted $\mb Q_\ell$-representation. Quotient stacks are
treated degree by degree through finite-dimensional approximations.
\end{proof}
Thus each Tate multiplicity space carries a finite-order rational
action. To compare its trace with a count over a finite field, we use
the explicit Jordan stratification below.

\subsection{The trace of the Jordan strata}

For a graded pure Tate
cohomology object $H$ with a finite-order action $u$, set
\[ \op{wt}_u(H;T)= \sum_r\Tr(u\mid H^{2r}(r))T^r. \]
For quotient stacks this is understood in $\mb Q((T^{-1}))$.

\begin{lemma}\label{lem:jordan-trace}
The stack $\mathfrak N_\alpha$ has pure cohomology of Tate type.
For $\alpha=\iota\beta$, we have
\begin{equation}\label{eq:jordan-trace} \op{wt}_a\bigl(H_c(\mathfrak N_\alpha);T\bigr) =J_{\beta,a}(T).\end{equation}
Consequently this trace, evaluated at $T=q$, equals
$|\mathfrak N_{\beta,a}(\mb F_q)|$ over every finite field.
\end{lemma}
\begin{proof}
We stratify by the nilpotent Jordan partitions at the unfolded vertices.
A stratum is the quotient of a vector space of arrow intertwiners by
the centralizer of the chosen nilpotents. The centralizer has a unipotent radical isomorphic as a variety to an
affine space, and reductive quotient a product of general linear
groups. The compactly supported cohomology of the stratum is therefore
pure of Tate type and concentrated in even degrees. Order the tuples
of Jordan partitions by the product of the dominance orders. Taking
successive unions of minimal $a$-orbits in this finite partially ordered
set gives an $a$-stable closed filtration, whose successive differences
are disjoint unions of Jordan strata. Localization has zero connecting
maps because they change cohomological degree by one. Thus the whole
stack is pure and even, and its trace is the sum of the traces of these
differences.

An orbit of Jordan partition strata of length greater than one
contributes zero to the trace. On an $a$-fixed tuple, the partitions are constant along each
vertex orbit. For an orbit of $d_i$ copies of $B\mr{GL}_m$, the trace of
cyclic permutation is obtained from the ordinary weight series of
$B\mr{GL}_m$ by replacing $T$ with $T^{d_i}$. This follows either
from a homogeneous tensor basis or from the symmetric algebra of the
Chern classes. There is no permutation sign, since these classes
are even. The contribution is
\[ T^{-d_i m^2}\prod_{j=1}^m(1-T^{-d_i j})^{-1}. \]
An automorphism of a complex affine space acts trivially on its
one-dimensional top compactly supported cohomology. Combining the
arrow, unipotent, and reductive factors gives exactly the summand
of \eqref{eq:jordan-coefficient}. This proves \eqref{eq:jordan-trace}
as an identity of Laurent expansions of rational functions. Proposition~\ref{prop:nil-jordan} gives the finite-field evaluation.
\end{proof}

\subsection{BPS cohomology and the PBW isomorphism}\label{subsec:bps-objects}
Let $\mc M_\alpha$ be the coarse moduli space of semisimple
$\Pi_{\Qbar}$-modules of dimension $\alpha$. We denote the
preprojective BPS sheaf, with trivial stability, by $\ms B_\alpha$.
It is the normalized mixed Hodge module defined by the central-loop
construction in \cite[Theorem~B]{DavisonStacks}.
We set
\[ \mc B_\alpha:=H_c(\mc M_\alpha,\ms B_\alpha) \]
and use this compactly supported normalization throughout.
The PBW isomorphism below reconstructs the cohomology of the
endomorphism stacks from these BPS cohomology groups.

Purity and the ordinary Kac-polynomial identity give
\begin{equation}\label{eq:primitive-normalization} \mc B_\alpha=\bigoplus_{j\geq0}B_{\alpha,j}\otimes\mb L^j, \qquad \sum_j\dim B_{\alpha,j}\,T^j=A_{\Qbar,\alpha}(T).\end{equation}
Here each $B_{\alpha,j}$ is a finite-dimensional rational vector space,
and only finitely many are nonzero; see
\cite[Sections~8.1--8.2, equations~(78)--(79)]{DavisonStacks}.
Equivalently, $B_{\alpha,j}$ is the rational vector space underlying
$H_c^{2j}(\mc M_\alpha,\ms B_\alpha)(j)$, whose Hodge type
is $(0,0)$. In particular $\mc B_{\epsilon_v}=\mb Q$.
The automorphism acts on the sheaf construction and sends $B_{\alpha,j}$ to
$B_{a\alpha,j}$. On an individual summand only the stabilizer of
$\alpha$ acts.

\begin{proposition}[Equivariant PBW for endomorphism stacks]
\label{prop:primitive-pbw}
Let $\Qbar$ be a finite acyclic quiver over $\mb C$, and let $a$
be a finite-order quiver automorphism.
For each dimension vector $\alpha$, let
$G_\alpha=\prod_{v\in\Qbar_0}\mr{GL}_{\alpha_v}(\mb C)$
act by simultaneous change of bases, and put
\begin{align*}
 \mathfrak E_\alpha
  &= [\{(M,f):M\in\Rep(\Qbar,\alpha),\ f\in\End_{\Qbar}(M)\}/G_\alpha],\\
 \mathfrak N_\alpha
  &= [\{(M,f):M\in\Rep(\Qbar,\alpha),\ f\in\End_{\Qbar}(M)
                                      \text{ nilpotent}\}/G_\alpha].
\end{align*}
Let $\mc B_\alpha=H_c(\mc M_\alpha,\ms B_\alpha)$ denote
ordinary preprojective BPS cohomology with trivial stability, in the
compactly supported normalization \eqref{eq:primitive-normalization}. Take rational
coefficients and set
\[ \mb L=\mb Q(-1)[-2],\qquad \mathsf H=\bigoplus_{r\ge0}\mb L^{-r}, \]
with trivial $a$-action on these scalar factors; $\mathsf H$ is the
graded dual of $H^*(B\mb G_m)$. Then there are
$a$-equivariant isomorphisms
\begin{align*}
 \bigoplus_{\alpha\ge0}H_c(\mathfrak E_\alpha,\mb Q)
  &\iso \op{Sym}\left(\bigoplus_{\alpha\ne0}
                            \mc B_\alpha\otimes\mathsf H\right),\\
 \bigoplus_{\alpha\ge0}H_c(\mathfrak N_\alpha,\mb Q)
  &\iso \op{Sym}\left(\bigoplus_{\alpha\ne0}
                   \mc B_\alpha\otimes\mb L^{-1}\mathsf H\right).
\end{align*}
These are isomorphisms of dimension-vector graded, cohomologically
graded mixed Hodge structures, completed degreewise in the negative
Tate direction. In each fixed dimension and cohomological degree the
spaces are finite-dimensional. The action sends degree $\alpha$ to
degree $a\alpha$, and the symmetric algebra is formed before taking
fixed degrees. The symmetric powers use the usual cohomological symmetry, which
has no odd signs here because the BPS cohomology groups are even. The
isomorphisms concern graded objects, not Hall algebra multiplications.
\end{proposition}
\begin{proof}
We apply the cohomological PBW theorem
\cite[Theorem~C]{DM}, in the form of
\cite[Section~3.4.2 and Theorem~3.11]{DavisonStacks}. We give the
compatibilities and normalizations needed for this application.

Adjoin a reverse arrow $b^*$ for each arrow $b:i\to j$ and a loop
$\omega_v$ at each vertex. For our left-to-right path convention the
potential on the resulting tripled quiver is
\[ \mc W=\sum_{b:i\to j}(\omega_jb^*b-\omega_i bb^*). \]
Extend $a$ by $a(b^*)=a(b)^*$ and $a(\omega_v)=\omega_{a(v)}$.
The potential is preserved, and its trace vanishes on its critical
locus by the Euler identity for a homogeneous cubic. We take the
stability for which all representations are semistable. It satisfies
the genericity hypothesis in the PBW theorem because the
antisymmetric Euler form of the tripled quiver is zero.

The condition that every $\omega_v$ be nilpotent defines a Serre
subcategory of representations of the \emph{tripled} quiver. To make
the restriction explicit, let $Z_\alpha$ be the closed locus in its
coarse moduli space where the characteristic polynomial of each
$\omega_v$ is $X^{\alpha_v}$. Characteristic polynomials multiply in
short exact sequences, so this condition is detected on
semisimplification. A direct sum belongs to $Z=\coprod_\alpha Z_\alpha$
if and only if both summands do. Base change for semisimplification
and for the finite direct-sum maps therefore restricts the relative
integrality isomorphism of
\cite[Theorem~3.9, equation~(22)]{DavisonStacks} to $Z$, with symmetric
powers of the restricted BPS sheaf on the right. This is the argument
in \cite[proof of Lemma~4.2, equation~(35)]{DavisonStacks}; dimensional
reduction on this locus is equation~(37) there. The automorphism $a$
preserves $Z$ and commutes with all these maps. We check the
equivariance of the PBW morphism below.

Let $\chi_{\Qbar}$ be the ordinary Euler form:
$\chi_{\Qbar}(\alpha,\eta)=\sum_v\alpha_v\eta_v-\sum_{b:i\to j}\alpha_i\eta_j$.
We choose the Hall sign twist $\psi=\chi_{\Qbar}$ modulo two. Indeed,
the Euler form of the tripled quiver satisfies
\[ \chi_{\Qbar^{\mr{tr}}}(\alpha,\eta) =-\sum_{b:i\to j}(\alpha_i\eta_j+\alpha_j\eta_i), \]
so its diagonal is even and its value modulo two equals
$\psi(\alpha,\eta)+\psi(\eta,\alpha)$. This is the symmetry correction
required in \cite[Theorem~C]{DM}; it is invariant under $a$.
For the scalar class in dimension $\alpha\ne0$, take
\[ z_\alpha=\frac1{|\alpha|} c_1\!\left(\bigotimes_v\det\mc V_v\right), \]
where $\mc V_v$ is the tautological vertex bundle. Its restriction
to the common scalar subgroup is the standard generator, and
$a(z_\alpha)=z_{a\alpha}$.

Let $\jmath$ be the BPS inclusion obtained by perverse truncation in
\cite[equation~(24)]{DavisonStacks}. The PBW map in Theorem~3.11 of that source is formed
from this inclusion, powers of $z_\alpha$, and symmetrized Hall
products. The isomorphisms induced by $a$ commute with
semisimplification and with the correspondences of short exact
sequences. Perverse truncation, the scalar action, and Hall pullback
and pushforward consequently commute with $a$. Thus
\[ a\bigl(z_\alpha^r\jmath(v)\bigr) =z_{a\alpha}^r\jmath(av),\qquad a(x*y)=(ax)*(ay). \]
For the PBW morphism this gives
$a\circ\Phi_{\mr{PBW}}=\Phi_{\mr{PBW}}\circ\op{Sym}(a\otimes\mr{id})$.
We apply this morphism in the category of cohomologically graded
mixed Hodge structures.

We next pass to compactly supported cohomology and compute the scalar
factors. For this comparison we use the \emph{half Tate object}
$\mb L^{1/2}$ in the monodromic category of the PBW theorem; its
square is $\mb L$. The sheaf identity in
\cite[Theorem~B and equation~(75)]{DavisonStacks} is stronger than
support on a single eigenvalue: under the map adjoining the scalar
loop $x\op{id}$, the tripled BPS sheaf is the pushforward
of
\[ \ms B_\alpha\boxtimes \bigl(\mb Q_{\mb A^1}\otimes\mb L^{-1/2}\bigr). \]
This identity is obtained from the zero-eigenvalue restriction and
translation of all vertex loops by $x$. The translation preserves the
trace function, since $d_M(x\op{id})=0$, and commutes with
$a$. It therefore preserves the vanishing-cycle construction and gives
the equivariant product with the scalar line, on which $a$ acts
trivially. Integrating
over $\mb A^1$ gives $\mc B_\alpha\otimes\mb L^{1/2}$;
restriction to $x=0$ gives
$\mc B_\alpha\otimes\mb L^{-1/2}$.
In the compactly supported PBW convention the common scalar factor is
$H_c(B\mb G_m)\otimes\mb L^{1/2}$. Since
\[ H_c(B\mb G_m)=\bigoplus_{r\ge1}\mb L^{-r} =\mb L^{-1}\mathsf H, \]
the two factors simplify to
$\mc B_\alpha\otimes\mathsf H \ \text{and}\ \mc B_\alpha\otimes\mb L^{-1}\mathsf H$.
The half twists cancel from both final expressions, which are
ordinary graded mixed Hodge structures.

Finally, dimensional reduction along the reverse arrows identifies
the normalized critical cohomology with the endomorphism-stack
cohomology. To check the normalization, put
$v=\dim\Rep(\Qbar,\alpha)$. The tripled
representation stack is smooth of dimension $2v$, so its normalized
constant object is $\mb Q\otimes\mb L^{-v}$. Reduction along the
$v$ reverse-arrow coordinates contributes $\mb L^v$ and cancels
that twist. It gives $H_c(\mathfrak E_\alpha)$, or
$H_c(\mathfrak N_\alpha)$ with the nilpotent-loop restriction; see
\cite[Theorem~3.7 and equations~(33)--(38)]{DavisonStacks}.
All the reduction maps commute with $a$, and graded duality preserves
this compatibility. Purity gives the even Tate grading. At fixed
$\alpha$ and cohomological degree only finitely many summands
contribute: there are at most $|\alpha|$ factors of positive dimension,
and each factor has a bounded-above Tate grading. Passing degreewise to the completion proves
the two isomorphisms.
\end{proof}

In the symmetric algebra, $a$ can permute factors of different
dimensions. Such factors may contribute a nonzero trace in their sum
degree even when none of their individual dimensions is fixed. We must
therefore form the symmetric algebra before restricting to invariant
dimensions.

\subsection{Norms, traces, and absolute indecomposability}

If $u$ is a power of $a$ and $u\alpha=\alpha$, write
\[ P_{\alpha,u}(T)=\sum_j\Tr(u\mid B_{\alpha,j})T^j\in\mb Z[T]. \]
Set this expression to zero if $u\alpha\ne\alpha$. Write $f_u$
for the inverse of repetition on the vertex orbits of $u$; in
particular $f=f_a$ is the inverse of $\iota$.

\begin{lemma}\label{lem:norm-log}
Put $\mc J_a(\mathbf x;T)=\sum_\beta J_{\beta,a}(T)\mathbf x^\beta$.
Then
\begin{equation}\label{eq:norm-cohomology-log} \log\mc J_a(\mathbf x;T)= \sum_{n\geq1}\frac{1}{n(T^n-1)} \sum_{\substack{\alpha\ne0\,:\ a^n\alpha=\alpha}} P_{\alpha,a^n}(T^n)\, \mathbf x^{f(\mc N_n^a\alpha)}.\end{equation}
Every coefficient is a finite sum. At a finite field the same identity
holds with $P_{\alpha,a^n}(q^n)$ replaced by $\mathsf A_{\alpha,n}(q)$.
\end{lemma}
\begin{proof}
For a finite-dimensional graded vector space, the trace identity for
symmetric powers is
\[ \sum_{r\ge0}\Tr(u\mid\op{Sym}^r V)z^r =\exp\left(\sum_{n\ge1}\frac{\Tr(u^n\mid V)}n z^n\right). \]
We apply this identity with variables recording dimension vectors and
cohomological degrees. If $\mathsf D_{\mathbf z}$ multiplies degree
$\alpha$ by $\mathbf z^\alpha$, then $(a\mathsf D_{\mathbf z})^n$ has a
trace on that degree only when $a^n\alpha=\alpha$. Its monomial is
$\mathbf z^{\mc N_n^a\alpha}$, not $\mathbf z^{n\alpha}$.
In the $n$-th logarithmic trace summand, the scalar factor
$\mb L^{-1}\mathsf H$ contributes
$\sum_{r\geq1}T^{-nr}=1/(T^n-1)$. The norm is invariant
when $a^n\alpha=\alpha$, so its monomial descends through $f$.
Proposition~\ref{prop:primitive-pbw} and
Lemma~\ref{lem:jordan-trace} now give
\eqref{eq:norm-cohomology-log}. We may first truncate the cohomological grading. Passing to the limit
then gives the identity as a Laurent series in $T^{-1}$, with
rational-function coefficients in each rank.

For the finite-field assertion, take the logarithm of
\eqref{eq:nil-product}. A geometric $F_a$-orbit of size $e$ is counted
$e$ times among the classes fixed by $F_a^n$ when $e\mid n$, and not
at all otherwise. For $n=e\ell$, its norm is $\ell$ times the
geometric dimension of the descended representation. Its contribution to the right-hand side is
$\frac{e}{e\ell(q^{e\ell}-1)}\mathbf x^{\ell\gamma} =\frac{\mathbf x^{\ell\gamma}} {\ell(q^{e\ell}-1)}$,
which is precisely its contribution to the logarithm of
\eqref{eq:nil-product}. Finally, a term contributing to rank $\beta$
satisfies $n|\alpha|=|\iota\beta|$, which proves coefficientwise
finiteness.
\end{proof}

\begin{example}\label{ex:vertex-orbit-trace}
Take $d$ isolated unfolded vertices $v_0,v_1,\dots,v_{d-1}$ with $a(v_s)=v_{s+1}\ (s\bmod d)$.
There are no quiver arrows; the permutation specifies the Galois
action. Their species is the one-vertex field $K_d$. For its finite-field analogue $\mb F_{q^d}$, the rank-one
stack of representations with a nilpotent endomorphism has cardinality $1/(q^d-1)$: the only module has
residue degree $d$, and its only nilpotent endomorphism is zero.
For $d>1$ its absolute count is zero. On the cohomological side
there are $d$ one-dimensional BPS spaces at the individual vertices,
none in an invariant dimension, and their length-$d$ cycle gives the term $1/(T^d-1)$ in their sum degree. At $d=1$ the scalar factors are
$1/(T-1)$ for nilpotent endomorphisms and $T/(T-1)$ for unrestricted
endomorphisms, fixing both normalizations.
\end{example}

\begin{theorem}[Absolute trace formula]\label{thm:allrank-absolute}
For every $\beta\ne0$ and every finite field, we have
\[ I_{\beta,1}(q)=P_{\iota\beta,a}(q) =\sum_j\Tr(a\mid B_{\iota\beta,j})q^j. \]
In particular, the absolute count is a polynomial in $\mb Z[T]$.
The formula also holds for the partially split species attached to
every power of $a$, with ranks measured on its vertex orbits.
\end{theorem}
\begin{proof}
We induct on $|\alpha|$, simultaneously over all finite fields and
all powers $u$ of $a$. Fix a $u$-invariant vector $\alpha$ and put
$\beta=f_u(\alpha)$. In the coefficient of $\mathbf x^\beta$ in the
$u$-version of \eqref{eq:norm-cohomology-log}, the $n=1$ term is
$P_{\alpha,u}(q)/(q-1)$; in the arithmetic version it is
$I^{(u)}_{\beta,1}(q)/(q-1)$. Here the superscript records the
species descended by $u\circ\mr{Frob}_q$.

Every remaining term involves $u^n$, $q^n$, and a dimension vector
$\eta$ with $n|\eta|=|\alpha|$, where $n\ge2$. These terms agree by
induction. The two left sides agree by Lemma~\ref{lem:jordan-trace}
and Proposition~\ref{prop:nil-jordan}; subtracting the remaining
terms proves equality of the $n=1$ terms. The induction starts at
total dimension one. Finally, a finite-order rational operator has a trace that is both rational and an algebraic integer.
Thus $P_{\alpha,u}$ has integer coefficients.
The comparison uses the rational-function identity of
Lemma~\ref{lem:jordan-trace} and its finite-field evaluation from
Proposition~\ref{prop:nil-jordan}, so it is valid in every finite
characteristic.
\end{proof}

\begin{corollary}\label{cor:residue-traces}
For all $\beta,e$, the refined count is the polynomial
\begin{equation}\label{eq:residue-traces} I_{\beta,e}(T)=\frac1e \sum_{\mc N_e^a\alpha=\iota\beta} \sum_{r\mid e}\mu_{\mr{ar}}(e/r) P_{\alpha,a^r}(T^r)\in\mb Q[T].\end{equation}
It is integer-valued on every finite-field cardinality. Its total
$I_{\beta,a}(T)=\sum_e I_{\beta,e}(T)$ is also in $\mb Q[T]$.
Remark~\ref{rem:nonintegral-count} gives a total count with
nonintegral coefficients.
\end{corollary}
\begin{proof}
Substitute Theorem~\ref{thm:allrank-absolute}, over each partial
splitting and extension field, into Lemma~\ref{lem:mobius}.
\end{proof}

\subsection{A signed trace polynomial}\label{subsec:signed-example}
\begin{example}\label{ex:arithmeticfamily}
Let $L=K_6=k((u))$, $u^6=t$, and take vertex fields $K_3,K_2$
with one arrow bimodule $L$, directed from the $K_3$-vertex to the
$K_2$-vertex:
\begin{equation} \label{eq:3,2}
\quiverdiagram{\xymatrix@C=3em{K_3\ar[r]^{K_6}&K_2}}
\qquad\qquad
\quiverdiagram{\xymatrix@C=2.5em@R=1.8em{
 &1\ar[dl]\ar[dr]&\\
 4&2\ar[l]\ar[r]&5\\
 &3\ar[ul]\ar[ur]&
}},\qquad a=(123)(45).
\end{equation}

The three sources $1,2,3$ lie above the $K_3$-vertex, and $4,5$
lie above the $K_2$-vertex. Then
\[ C=\begin{pmatrix}2&-2\\-3&2\end{pmatrix},\qquad D=\op{diag}(3,2). \]
At rank $(1,1)$, a nonzero arrow is
$x\mapsto\Tr_{L/K_2}(cx)$, $c\in L^\times$, and its class is the coset
$c\in L^\times/(K_3^\times K_2^\times)$.
Each such representation is absolutely indecomposable and has
endomorphism field $K$. More geometrically, let
$W=\Res_{L/K}\mb G_m\subset V_{(1,1)}$. This is the open
subset on which the field norm of $c$ is nonzero. After extending
scalars, all six arrow coordinates are nonzero. Since the quiver is
connected and all vertex spaces are one-dimensional, its endomorphisms
are scalar. The family is smooth, geometrically integral and invariant
under change of bases, and it contains the point $c=1$. Its dimension is $[L:K]=6$,
whereas $\dim G_{(1,1)}=\lVert(1,1)\rVert_d^2=3+2=5$. Hence
\[ q_d(1,1)=3+2-6=-1,\qquad \nu(W)=6-5+1=2=1-q_d(1,1). \]
Thus this example also realizes the parameter dimension in
Theorem~\ref{thm:absolute-root}.

The elements $1+zu$, $z\in k$, represent distinct cosets. Indeed,
a ratio $(1+zu)/(1+z'u)$ belonging to $K_3^\times K_2^\times$
can be written as a product of principal units from the two subfields
after cancelling valuations and constant terms. The logarithm of such
a product has no linear $u$-term, whereas the logarithm of the ratio
has linear term $(z-z')u$. Hence $z=z'$. Over $\mb F_q$, the same quotient of multiplicative groups gives
\[ I_{(1,1),a}(q)=I_{(1,1),1}(q) =\frac{(q^6-1)(q-1)}{(q^3-1)(q^2-1)}=q^2-q+1. \]

The unfolding is $\mathsf K_{3,2}$ and $\alpha=(1,1,1;1,1)$. Connected
nonzero supports give twelve trees, six five-edge subgraphs and the
full graph, so
\[ A_{\Qbar,\alpha}(T)=12+6(T-1)+(T-1)^2=T^2+4T+7. \]
Both the ordinary Kac polynomial and the species trace polynomial
have degree two, equal to $\nu(W)$, although their coefficients differ.
The trace polynomials for $a=(123)(45)$ and its powers are
\[
\begin{array}{c|cccc}
 r&0&1,5&2,4&3\\\hline
 P_{\alpha,a^r}(T)&T^2+4T+7&T^2-T+1&T^2+T+1&T^2+2T+1.
\end{array}
\]
For $r=2$, the two edge orbits are three-edge stars; for $r=3$ there
are three two-edge orbits. Connectivity requires all of them, and
the scalar quotients give $(q^3-1)/(q-1)$ and $(q+1)^2$, respectively.
Let $\varepsilon$ be the sign representation of the cyclic group of
order six, and $\rho_3$ the two-dimensional rational representation whose eigenvalues
are the primitive cube roots of unity. The table and Theorem~\ref{thm:allrank-absolute} give
\[ B_{\alpha,2}=\mb Q,\qquad B_{\alpha,1}\iso\mb Q\oplus\varepsilon\oplus\rho_3,\qquad B_{\alpha,0}\iso\mb Q\oplus\mb Q[\mb Z/6\mb Z]. \]
Their dimension polynomial is $T^2+4T+7$, their trace at $a$ is
$T^2-T+1$, and the polynomial of invariant dimensions is $T^2+T+2$.
Thus the trace at $a$ differs both from the dimension and from the
dimension of the invariant subspace.
\end{example}

\subsection{A recursion using partial splittings}

The proof also gives a finite recursion for the BPS trace polynomials
without requiring an explicit quiver variety. For an automorphism $u$ write
$f_u$ for the inverse of the repetition map on its vertex orbits and
$\mc J_u$ for its Jordan series. The notation $[\mathbf x^\gamma]$
means the coefficient of that monomial. At a $u$-invariant vector
$\alpha\ne0$, isolate the $n=1$ term to obtain
\begin{equation}        \label{eq:primitive-recursion}
 P_{\alpha,u}(T)=(T-1)\Bigg([\mathbf x^{f_u(\alpha)}] \log\mc J_u(\mathbf x;T)\
 -\ \sum_{\substack{n\geq2\,,\ \eta\ne0\\
                         \mc N_n^u\eta=\alpha}}
       \frac{P_{\eta,u^n}(T^n)}{n(T^n-1)}\Bigg).
\end{equation}
Each remaining term has smaller total unfolded dimension, so
\eqref{eq:primitive-recursion} is a finite recursion at any prescribed
rank bound. Theorem~\ref{thm:allrank-absolute} shows that the result
lies in $\mb Z[T]$.

The two modulations of Example~\ref{ex:two-modulations} have equal
rank-$\delta$ contributions to $\mc J_a$,
$\frac{T^2+1}{T^2-1}=\frac{T-1}{T-1}+\frac{2}{T^2-1}$, although their
indecomposables have different residue degrees. Thus $\mc J_a$
alone does not determine the refined counts. Their $a^2$-splittings
are a disjoint union of two Kronecker quivers and a connected square,
respectively; the recursion distinguishes them.

\subsection{The degree-zero action and twining}

We use the rational Chevalley--Serre forms
$\mathfrak n_{\mb Q}^{\pm}(\widetilde C)$ of the positive and
negative halves of the Kac--Moody algebra. Their complexifications
are the usual complex Lie algebras. The diagram action sends the
positive generator $E_v$ to $E_{a(v)}$ and acts in the same way on the
negative generators. All traces in this subsection are on these
finite-dimensional rational root spaces.

\begin{lemma}\label{lem:bps-zero-action}
In the normalization \eqref{eq:primitive-normalization}, we have
\begin{equation}\label{eq:bps-zero-action} B_{\alpha,0}\iso \bigl(\mathfrak n_{\mb Q}^+(\widetilde C)_\alpha\bigr)^\vee\end{equation}
compatibly with the quiver automorphism. Consequently
\[ P_{\alpha,a^n}(0) =\Tr(a^n\mid\mathfrak n_{\mb Q}^+(\widetilde C)_\alpha) \quad\text{when }a^n\alpha=\alpha. \]
\end{lemma}
\begin{proof}
We first match the two conventions for BPS cohomology. The
preprojective BPS sheaf in
\cite[Theorem/Definition~4.1, equations~(40)--(41)]{DavisonBPS}
is characterized by the scalar-loop product identity. Restricting
to the zero scalar eigenvalue, with the half Tate correction,
identifies it with $\ms B_\alpha$ above by
\cite[Theorem~B, equation~(7)]{DavisonStacks}.
The latter theorem also gives Verdier self-duality; its proof in
\cite[Section~4.2, equation~(43)]{DavisonStacks} uses the natural
dualities of the normalized intersection complex and vanishing
cycles, together with the scalar-line product decomposition.
The automorphism $a$ preserves the potential, commutes with
scalar-loop translation, and preserves the zero scalar eigenvalue.
These duality and product identifications are therefore compatible
with $a$. Verdier duality gives the equivariant isomorphism
\[ \mc B_\alpha^\vee
   =H_c(\mc M_\alpha,\ms B_\alpha)^\vee
   \iso H(\mc M_\alpha,\ms B_\alpha). \]
The last group is the dimension-$\alpha$ summand of the BPS Lie
algebra in \cite[Definition~4.2]{DavisonBPS}, with the Hall sign
convention $\chi_{\Qbar}$ used in
Proposition~\ref{prop:primitive-pbw}. In cohomological degree zero
it is $B_{\alpha,0}^{\vee}$.

Davison's Theorem~6.6 gives an isomorphism over $\mb Q$
\[ \mathfrak n_{\mb Q}^-(\widetilde C) \longrightarrow H^0(\mathfrak g_{\Pi_{\Qbar}}) \]
by sending the Chevalley generators to the canonical vertex classes.
Its loopless subquiver is all of $\Qbar$, so the statement includes
Dynkin quivers. The automorphism $a$ permutes those classes and
preserves the Hall product with the chosen sign convention. The two
composites of this isomorphism with $a$ agree on every generator,
and hence agree on the entire Lie algebra. Identifying the negative
and positive Chevalley--Serre presentations and taking the graded
dual proves \eqref{eq:bps-zero-action}.

For a finite-order rational representation the traces of an element
and its inverse are equal. Indeed, its eigenvalues are roots of unity
and its rational character is fixed by complex conjugation. This
removes the dual in the trace identity.
\end{proof}

We record the root-lattice convention for the twining theorem. For an
admissible automorphism of a symmetric Cartan matrix, its orbit
algebra, in the standard convention used here, has entries
$C_{ij}=\sum_{w\in O_j}\widetilde C_{v w}$ for $v\in O_i$.
It has symmetrizer $\op{diag}(|O_i|)$, and its simple root
pulls back to $\sum_{v\in O_i}\epsilon_v=\iota\epsilon_i$.
Indeed the product of commuting reflections over $O_i$ restricts to
$s_i$ on the invariant lattice. The Cartan convention in
\cite[equation~(3.1)]{FSS} is transposed: their Serre exponent is
$1-A^{j,i}$. Transposing their orbit-matrix formula, equation~(2.10),
gives the entries $C_{ij}$ above. Their linking factors are all one
because $a$ is admissible; equation~(5.12) then identifies the pullback
of a simple root with the sum over its orbit. We use the orbit Lie algebra in the twining theorem, not
the fixed-point subalgebra.

For an $a$-equivariant multigraded vector space
$V=\bigoplus_\alpha V_\alpha$, define its \emph{twining character} by
\[ \op{ch}_a(V)=\sum_{\beta\in\mb N^I}
   \Tr\!\left(a\mid V_{\iota\beta}\right)\mathbf x^\beta. \] The twining theorem for a Verma module of
highest weight zero \cite[Theorem~1(a)]{FSS}, together with the ordinary
PBW theorem, says
\begin{equation}\label{eq:twining-denominator} \op{ch}_a\!\left( \op{Sym}\mathfrak n_{\mb Q}^+(\widetilde C)\right) =\prod_{\beta\in\Delta^+(C)} (1-\mathbf x^\beta)^{-\dim\mathfrak g(C)_\beta}.\end{equation}
If $a$ permutes several connected components,
the trace of the cyclically permuted tensor factors reduces this
identity to one component with its return automorphism. Products over
orbits of components then give \eqref{eq:twining-denominator} in general.

\begin{theorem}[Folded constant term]\label{thm:folded-constant}
For every nonzero rank $\beta$, we have
\begin{equation}\label{eq:folded-constant} I_{\beta,a}(0)=\dim\mathfrak g(C)_\beta,\end{equation}
where the right side is zero when $\beta$ is not a root. This is the
constant term of the count of \emph{all} indecomposables, not only the
absolute count. It is independent of the integral symmetrizer and the
ramified modulation realizing $C$. Equivalently,
\begin{equation}\label{eq:Jordan-denominator-zero} \mc J_a(\mathbf x;0) =\prod_{\beta\in\Delta^+(C)} (1-\mathbf x^\beta)^{\dim\mathfrak g(C)_\beta}.\end{equation}
Every coefficient of the rational Jordan series is regular at $T=0$
after simplification, so this specialization is well-defined.
\end{theorem}
\begin{proof}
By Lemma~\ref{lem:norm-log}, a coefficient of
$\log\mc J_a(\mathbf x;T)$ is a finite sum of polynomials divided
by $n(T^n-1)$. It is regular at zero, and so is each coefficient of its
formal exponential. Using Lemma~\ref{lem:bps-zero-action}, its value at
zero is
\[ -\sum_{n\geq1}\frac1n \sum_{a^n\alpha=\alpha} \Tr(a^n\mid\mathfrak n_{\mb Q}^+(\widetilde C)_\alpha) \mathbf x^{f(\mc N_n^a\alpha)}. \]
This is the negative logarithm of the left side of
\eqref{eq:twining-denominator}, by the trace identity of Lemma~\ref{lem:norm-log}.
The twining theorem proves \eqref{eq:Jordan-denominator-zero}.

On the other hand, \eqref{eq:nil-product} and
Corollary~\ref{cor:residue-traces} give a rational-function identity
coefficient by coefficient,
\[ \log\mc J_a(\mathbf x;T) =\sum_{\gamma\ne0}\sum_e I_{\gamma,e}(T) \sum_{n\geq1}\frac{\mathbf x^{n\gamma}} {n(T^{en}-1)}. \]
This identity holds over $\mb Q(T)$ because it holds at every finite
field. At zero it becomes
\[ -\sum_{\gamma\ne0} I_{\gamma,a}(0) \sum_{n\geq1}\frac{\mathbf x^{n\gamma}}n. \]
The exponents in a formal Euler product are unique, as can be seen
inductively in total rank, or by M\"obius inversion along each ray.
Comparing with \eqref{eq:Jordan-denominator-zero} proves
\eqref{eq:folded-constant}.
\end{proof}

The twisted degree-$(2,2)$ affine species illustrates the need for
the total count: its absolute constant is $-1$, while the degree-two
contribution has constant $2$, giving the root multiplicity $1$.

\begin{corollary}\label{cor:affine-lie-multiplicity}
In the notation of Corollary~\ref{cor:affine-counts}, we have
\[ \dim\mathfrak g(C)_{m\delta} =1+\sum_{\OO:\ h_{\OO}\mid m} (g_{\OO}-1). \]
Although the separate tube data can depend on the modulation, this
combination does not.
\end{corollary}
\begin{proof}
Combine \eqref{eq:affine-constant-count} with
Theorem~\ref{thm:folded-constant}.
\end{proof}

\subsection{A smooth realization}
When $g_d(\beta)=1$, the absolute counting polynomial
$I_{\beta,1}(T)$ admits the following realization on an ordinary
smooth moduli space. This is the equivariant form of
the construction in \cite{CBVdB}.

\begin{corollary}\label{cor:smooth-trace}
Suppose $k=\mb C$, $\beta\in\Delta^+(C)$, and $g_d(\beta)=1$.
Put $\alpha=\iota\beta$ and $p=1-q_d(\beta)$, and choose
$\lambda\in\mb Z^I$ such that $\bar\lambda=\iota\lambda$ is
$\alpha$-generic. Let $\Pi^0_{\Qbar}$ denote the ordinary preprojective algebra of
$\Qbar$, and set
\[ X^{\mr{st}}= \Rep(\Pi^0_{\Qbar},\alpha)^{\bar\lambda\text{-}\mr{st}} /(G_\alpha/\mb G_m), \]
where stability means $\bar\lambda\cdot\dimv U>0$ for every nonzero
proper subrepresentation $U$. Then $X^{\mr{st}}$ is smooth of
dimension $2p$, and
\begin{equation}\label{eq:smooth-trace} I_{\beta,a}(T)=I_{\beta,1}(T) =\sum_{j=0}^{p} \Tr\!\left(a\mid H^{2p-2j}(X^{\mr{st}},\mb Q)\right)T^j.\end{equation}
Under the equivariant comparison with the deformed quotient
$X_K^\lambda=\Rep(\Pi_A^\lambda,\beta)/(G_\beta/\mb G_m)$ and
comparison with $\mb Q_\ell$-cohomology, this finite action is its
inertia action, up to inversion of the chosen generator.
\end{corollary}
\begin{proof}
Propositions~\ref{prop:root}, \ref{prop:genericity}, and
\ref{prop:divisibility} give the parameter and an absolute model,
so $\alpha$ is an indivisible unfolded root. If $\alpha$ is a vertex
simple, both moduli spaces are points. Otherwise the smooth family
over the parameter line
$\mb A^1\bar\lambda$ in \cite[Section~2.1]{CBVdB} has zero
fibre $X^{\mr{st}}$. Its fibre at $\bar\lambda$ is the deformed
quotient $X=\Rep(\Pi^{\bar\lambda}_{\Qbar},\alpha)/(G_\alpha/\mb G_m)$.
The automorphism $a$ acts on this family and fixes the parameter line
pointwise. The cohomological comparison in
\cite[Proposition~2.3.1 and Corollary~2.3.2]{CBVdB} is induced by the adjunction morphism for the projection from the
parameter line to the arithmetic base. It is therefore natural under
automorphisms of the family, and commutes with $a$ and, in good
characteristic, with Frobenius.

For the $a\circ\mr{Frob}_q$-form of the deformed quotient, the
lifting criterion and scalar stabilizer give
$|X^{(a)}(\mb F_q)|=q^p I_{\beta,1}(q)$; this is the calculation
of \cite[Proposition~2.2.1]{CBVdB}, using the field-trace pairing.
Lang's theorem applies to the connected group $G_\beta/\mb G_m$.
Purity and the untwisted point counts over finite extensions in
\cite[Section~2.4]{CBVdB} give only even cohomology, with Frobenius
eigenvalues $q^r$ in degree $2r$. Since Frobenius commutes with $a$, a possible unipotent part makes no
contribution to the trace on any $a$-eigenspace. The equivariant
comparison and the fixed-point formula therefore give
\eqref{eq:smooth-trace} in
sufficiently large characteristic; traces of $a$ and $a^{-1}$ agree
over $\mb Q$. Both sides are polynomials, by
Theorem~\ref{thm:allrank-absolute}, so the identity holds in
$\mb Z[T]$. Finally $X_K^\lambda$ becomes constant over $K_N$,
with descent action $a$, giving the inertia assertion.
\end{proof}

\section{Counting polynomials and arithmetic isomorphism classes}\label{sec:polynomial-strata}
We compare the degrees and leading coefficients of the counting
polynomials with the geometric classification of Section~\ref{sec:classification}.

\subsection{Degrees and cardinalities}

\begin{lemma}\label{lem:leading-action}
If $\alpha$ is a positive root fixed by $u$, the polynomial
$P_{\alpha,u}(T)$ is monic of degree $p(\alpha)$. At a nonroot
it vanishes.
\end{lemma}
\begin{proof}
The ordinary Kac polynomial is monic of degree $p(\alpha)$, so its
top multiplicity space is a rational line. A finite-order operator
on that line is $1$ or $-1$. By Theorem~\ref{thm:allrank-absolute},
$P_{\alpha,u}(q)$ is nonnegative for arbitrarily large prime powers;
its leading coefficient therefore cannot be $-1$. At a nonroot,
\eqref{eq:primitive-normalization} gives $\mc B_\alpha=0$.
\end{proof}

\begin{corollary}\label{cor:allrank-polynomial-strata}
For every nonzero rank $\beta$ and residue degree $e$, one has
\begin{align}
 I_{\beta,e}(T)
 &=|\ms R^{\mr{re,exact}}_{\beta,e}|\ +\ \sum_{[\alpha]\in\ms R^{\mr{im}}_{\beta,e}}
   \frac{s(\alpha)}e
   \sum_{\substack{r\mid e\\s(\alpha)\mid r}}
      \mu_{\rm ar}(e/r)P_{\alpha,a^r}(T^r).
       \label{eq:allrank-polynomial-orbits}
\end{align}
Write $\op{lc}$ for the leading coefficient. If
$\ms R^{\mr{im}}_{\beta,e}\ne\varnothing$, then
\begin{equation}\label{eq:allrank-degree-leading}
 \deg I_{\beta,e}=e\max_{[\alpha]\in\ms R^{\mr{im}}_{\beta,e}}
                          p(\alpha),\qquad
 \op{lc} I_{\beta,e}
   =\sum_{\substack{[\alpha]\in\ms R^{\mr{im}}_{\beta,e}\\
              e p(\alpha)=\deg I_{\beta,e}}}\frac{s(\alpha)}e>0.
\end{equation}
Consequently $\Ind_{\beta,e}(\Lambda)$ is empty, finite nonempty,
or infinite exactly when $I_{\beta,e}(T)$ is zero, a nonzero
constant, or nonconstant, respectively. In the finite case the
constant is its exact cardinality.
\end{corollary}
\begin{proof}
We group the terms of the residue inversion
\eqref{eq:residue-traces} into $a$-orbits. An orbit of length $s$
contributes $s/e$ times the divisor sum in
\eqref{eq:allrank-polynomial-orbits}; terms not divisible by $s$
are zero because their automorphism does not fix the dimension.
For a real root all remaining trace polynomials equal one and
\[
\sum_{\substack{r\mid e\\s\mid r}}\mu_{\rm ar}(e/r)
                       =\begin{cases}1,&s=e,\\0,&s<e.\end{cases}
\]
This gives the first term. For an imaginary root $p(\alpha)>0$,
the term $r=e$ has degree $e p(\alpha)$ and leading coefficient
one before multiplication by $s/e$, by Lemma~\ref{lem:leading-action}. Every proper divisor gives
strictly smaller degree. Therefore the top coefficients of
all imaginary orbit contributions are positive and cannot cancel.
This proves \eqref{eq:allrank-degree-leading}. The trichotomy
now follows from Theorem~\ref{thm:allrank-strata}.
\end{proof}

\begin{corollary}\label{cor:primitive-specialization}
If $\beta$ is primitive, every representative $\alpha$ of an orbit
in $\ms R_{\beta,e}$ is primitive and has dimension period $e$. Thus
\[ I_{\beta,e}(T)= \sum_{[\alpha]\in\ms R_{\beta,e}} P_{\alpha,a^e}(T^e)\in\mb Z[T^e]. \]
\end{corollary}
\begin{proof}
A coordinate divisor of $\alpha$ divides $\beta$. A shorter
period $s<e$ would give
$\iota\beta=(e/s)\mc N_s^a\alpha$, contradicting
primitivity. All terms indexed by proper divisors vanish and $s/e=1$.
\end{proof}

At nonprimitive rank, the leading coefficient in
\eqref{eq:allrank-degree-leading} can be fractional, and the
smaller exponents need not be divisible by $e$.

For a root orbit of dimension period $s$ contributing at
residue degree $e=sm$, the factor $s/e=1/m$ in
\eqref{eq:allrank-degree-leading} counts Galois orbits of size $m$
over the return field rather than their individual members.
This explains the rational leading coefficients.

\subsection{Two ranks in an indefinite example}\label{subsec:indefinite-ranks}

We compare the primitive rank $(2,3)$ with its double $(4,6)$
for the degree-$(3,2)$ species of Example~\ref{ex:arithmeticfamily}.
We use the same vertex labels \eqref{eq:3,2}.

\begin{example}\label{ex:indefinite-four-degrees}
Take $\beta=(2,3)$. Although this rank vector is primitive,
$g_d(\beta)=6$, so primitivity does not force residue degree one.
Applying \eqref{eq:jordan-coefficient} and
\eqref{eq:primitive-recursion} gives the following polynomials.
\[
\begin{array}{c|l}
 e&I_{(2,3),e}(T)\\ \hline
 1&T^7-T^6+T^5+T-1\\
 2&T^4+T^2+1\\
 3&T^3+1\\
 6&1
\end{array}
\]
Corollary~\ref{cor:allrank-polynomial-strata} shows that the sets
of residue degrees one, two and three have cardinality $|K|$, whereas
there is exactly one class of residue degree six. This last class
comes from the orbit of the six real roots supported on single
arrows.
The total polynomial is
\[ T^7-T^6+T^5+T^4+T^3+T^2+T+2. \]
Its constant term is two, whereas the constant term of
$I_{(2,3),1}$ is $-1$.
The root space of degree $(2,3)$ for
$C=\left(\begin{smallmatrix}2&-2\\-3&2\end{smallmatrix}\right)$
also has dimension two: the free Lie algebra has dimension
$\binom{5}{2}/5=2$ in this multidegree, and the Serre relations start
at degrees $(3,1)$ and $(1,4)$, neither of which is bounded by
$(2,3)$. The equality of these constant terms is a special case of
Theorem~\ref{thm:folded-constant}.
\end{example}

\begin{example}\label{ex:doubled-indefinite}
Now take $\beta=(4,6)$.
Here $g_d(\beta)=12$, so only divisors of $12$ can occur as residue
degrees. The table gives the imaginary root orbits contributing to
each such degree. We list the quadratic value and dimension period
once for each orbit.
\[
\begin{array}{c|c|c|c|c}
 e&\bigl(q_{\Qbar}(\alpha),s(\alpha)\bigr)&
 \Ind_{(4,6),e}&\deg I_{(4,6),e}&\op{lc} I_{(4,6),e}\\\hline
 1&(-24,1)&\text{infinite}&25&1\\
 2&(-6,1),\ (-4,2)&\text{infinite}&14&1/2\\
 3&(-2,3),\ (0,3)&\text{infinite}&9&1\\
 4&(-1,2)&\text{infinite}&8&1/2\\
 6&(0,3)&\text{infinite}&6&1/2\\
 12&\varnothing&\text{empty}&\text{--}&\text{--}
\end{array}
\]
The absence of an imaginary root orbit for $e=12$ does not
by itself explain emptiness: we must also examine the real roots.
Let $\alpha_0=(1,0,0;1,0)$. Its orbit consists of the six roots
supported on single arrows, and
\[ s(\alpha_0)=6,\qquad \mc N_6^a\alpha_0=\iota(2,3),\qquad \mc N_{12}^a\alpha_0=\iota(4,6). \]
These six roots are all the solutions of the degree-twelve norm
equation. The unique indecomposable of each real-root dimension has
isomorphism period six, so none can have period twelve. Thus the
orbit that gives one class at rank $(2,3)$ and residue degree six
gives no class at rank $(4,6)$ and residue degree twelve.
By contrast, imaginary-root families can have an isomorphism period
larger than their dimension period. This occurs in the degree-two,
degree-four and degree-six contributions in the table.
In degree one the root $\alpha=(4,4,4;6,6)$ has
$q_{\Qbar}(\alpha)=-24$. It illustrates the construction at a
nonprimitive rank with $q_d(\beta)<0$.

For residue degree three, the norm equation has two root orbits,
represented by $(2,1,1;2,2)$ and $(2,2,0;2,2)$, each of period three.
The partial-species Jordan formula gives
\[ P_{(2,1,1;2,2),a^3}(T)=T^3+2T^2+4T+3, \qquad P_{(2,2,0;2,2),a^3}(T)=T+1.\]
There are no terms from proper divisors of three, so
$I_{(4,6),3}(T)=T^9+2T^6+5T^3+4$.
Three further residue degrees give simple polynomials:
\[ I_{(4,6),4}(T)=\frac{T^8-T^2}{2},\qquad I_{(4,6),6}(T)=\frac{T^6-T^3}{2}+1,\qquad I_{(4,6),12}(T)=0. \]
The total constant is $I_{(4,6),a}(0)=6$ by
Theorem~\ref{thm:folded-constant}. The Weyl denominator has truncation
$1-x-y+x^3y+xy^4$ through bidegree $(4,6)$, giving an independent check:
\[
[x^4y^6]\log(1-x-y+x^3y+xy^4)
=-7=-\dim\mathfrak g(C)_{(4,6)}-\tfrac12\dim\mathfrak g(C)_{(2,3)}.
\]
Since the latter multiplicity is $2$, the former is $6$.

The degree-four polynomial also has a direct geometric
interpretation. For $e=4$, the norm equation has one positive-root
orbit, represented by $\alpha=(1,1,1;1,2)$, whose dimension period is
$s(\alpha)=2$. Put $F=K_2$ and $L=K_6$, so $[L:F]=3$.
For the partially split species over $F$, a representation of this
dimension consists of $F$-linear maps
\[
  \quiverdiagram{\xymatrix@C=3.8em{
     F&L\ar[l]_{f}\ar[r]^{g}&F^2
  }}.
\]
It is absolutely indecomposable exactly when $f\ne0$ and $g$ is
surjective. Necessity follows because otherwise a target simple
splits off. Conversely, after extending to an algebraic closure,
$f\ne0$ has a nonzero component at every one of the three source
vertices. Hence an endomorphism has a common scalar on those vertices
and the first target; surjectivity of $g$ forces the same scalar on the
second target.

By the nondegenerate trace pairing and a source change of basis, we
may normalize $f(x)=\Tr_{L/F}(x)$. We can then write
\[
  g(x)=\bigl(\Tr_{L/F}(u_1x),\Tr_{L/F}(u_2x)\bigr),
  \qquad u_1,u_2\in L.
\]
Surjectivity means that $u_1,u_2$ are $F$-linearly independent.
Changing bases at the second target changes only their basis, while
the remaining source changes preserving the normalized first map up
to a target scalar lie in $F^\times$. Thus the isomorphism classes
are parametrized by
\[
  U=\langle u_1,u_2\rangle_F
    \in\op{Gr}_F(2,L)\cong\mb P^2_F.
\]
The same normal form holds over $K_4/F$, since this quadratic
extension is linearly disjoint from the cubic extension $L/F$. By Galois orbit descent, rank-$(4,6)$
indecomposables of residue degree four correspond to degree-two
closed points of this projective plane. The corresponding
finite-field count is therefore
\begin{equation*}
  I_{(4,6),4}(q)
   =\frac{|\mb P^2(\mb F_{q^4})|
           -|\mb P^2(\mb F_{q^2})|}{2}
   =\frac{q^8-q^2}{2}.
\end{equation*}
This also explains why the rank-$(2,3)$ polynomial at residue degree
two is $I_{(2,3),2}(q)=|\mb P^2(\mb F_{q^2})|=q^4+q^2+1$:
doubling both the rank and residue degree replaces rational points of
the same parameter space by its degree-two closed points.
\end{example}

\section{Integral orders and the special fibre}\label{sec:orders}
We compare lattices over integral orders with their generic fibres
and with locally free modules over the special-fibre GLS algebra.

\subsection{Decompositions}\label{subsec:lattice-splitting}

Put $\OO=k[[t]]$, $\OO_i=k[[t^{1/d_i}]]$, and
$\RR=\prod_i\OO_i$. A \emph{full order} $\ms A\subset\Lambda$ is an $\OO$-subalgebra
finite over $\OO$ with $K\otimes_\OO\ms A=\Lambda$.
We call $\ms A/t\ms A$ its \emph{special fibre} and
$K\otimes_\OO\ms A$ its \emph{generic fibre}. We assume
that it contains $\RR$. A \emph{lattice} is a right $\ms A$-module that is finite
and torsion-free over $\OO$. For each $i$, the module $Le_i$ is free over $\OO_i$, since $\OO_i$
is a discrete valuation ring.
Every generic representation has a lattice: if $L_0$ is a full
$\OO$-lattice in it, then $L_0\ms A$ is one. Clearing denominators gives
$K\otimes_\OO\Hom_{\ms A}(L,L') \iso\Hom_\Lambda(L_K,L'_K)$. 
Thus generic scalar extension gives an equivalence from the
\emph{isogeny category} of lattices (with morphism spaces tensored
with $K$) to $\modfd(\Lambda)$. An indecomposable $\ms A$-lattice
may nevertheless have a decomposable generic fibre.

\begin{example}[Localization need not preserve indecomposability]
\label{ex:localization}
For the Kronecker order $\OO Q$, take
\[ \ms M:\quad\quiverdiagram{\xymatrix@C=4em{\OO^2\ar@/^0.65pc/[r]^{I_2}\ar@/_0.65pc/[r]_{U_t}&\OO^2}},\qquad U_t=\begin{pmatrix}0&1\\0&t\end{pmatrix}. \]

Its endomorphisms are the matrices
$\bigl(\begin{smallmatrix}x&y\\0&x+ty\end{smallmatrix}\bigr)$, so
$\End_{\OO Q}(\ms M)\iso\OO[z]/(z^2-tz)$.
Direct reduction gives $\End_{kQ}(\ms M/t\ms M)\iso k[z]/(z^2)$.
The special fibre, and therefore the lattice, is indecomposable. Its generic fibre is the sum of two rank-$(1,1)$ Kronecker
representations:
\[ \ms M_K\iso
\left(\quiverdiagram{\xymatrix@C=2.6em{K\ar@/^0.55pc/[r]^{1}\ar@/_0.55pc/[r]_{0}&K}}\right)
\oplus
\left(\quiverdiagram{\xymatrix@C=2.6em{K\ar@/^0.55pc/[r]^{1}\ar@/_0.55pc/[r]_{t}&K}}\right). \]

The projector onto the second summand is
$U_t/t=\begin{pmatrix}0&t^{-1}\\0&1\end{pmatrix}$.
This disproves \cite[Lemma~5.17]{GLSrigid} as stated, since the
lemma imposes no rigidity hypothesis. Its intersection argument fails
to preserve a complementary sum: inside
$E=\{(x,y)\in\OO^2:x\equiv y\pmod t\}$, the intersections
with the two generic coordinate axes sum to $t\OO^2\subsetneq E$.
The example is nonrigid; it does not by itself contradict the rigid
main conclusions of that paper.
\end{example}

For the standard species of Proposition~\ref{prop:standard}, the order
of \cite[Sections~4.1--4.2]{GLSrigid} is
\[ \ms H=T_{\RR}\!\left(\bigoplus_{i\to j} (\OO_i\otimes_{\OO_{ij}}\OO_j)^{p_{ij}}\right), \qquad \OO_{ij}=k[[t^{1/h_{ij}}]]. \]
The arrow lattice $\OO_i\otimes_{\OO_{ij}}\OO_j$ need not be
the integral closure of $\OO$ in $K_{\ell_{ij}}$. The special fibre of $\ms H$ is
$H=H(C,D,\Omega)$, where $\Omega$ is the acyclic orientation.
Its vertex algebra is $H_i=k[\varepsilon_i]/(\varepsilon_i^{d_i})$.
An $H$-module $X$ is locally free when $Xe_i$ is free over $H_i$;
its rank at $i$ is that free rank. Every $\ms H$-lattice is free over each vertex order $\OO_i$.
By \cite[Lemma~5.4]{GLSrigid}, every locally free $H$-module admits
an $\ms H$-lattice lift, obtained by lifting the arrow maps.

\begin{corollary}\label{cor:nonroot-lifts}
Let $X$ be a nonzero locally free $H$-module whose rank is not a
root of $C$. Every lattice lift of $X$ has decomposable generic
fibre. If $X$ is indecomposable, every such lattice is indecomposable,
and every nontrivial generic splitting idempotent is nonintegral.
\end{corollary}
\begin{proof}
The generic fibre has the same rank, so Proposition~\ref{prop:root} forces
its decomposition. A decomposition of a lattice into two nonzero
summands would reduce to such a decomposition of its special fibre,
by Nakayama's lemma. A nontrivial integral idempotent would produce
exactly that decomposition.
\end{proof}

In particular, the indecomposable locally free modules of nonroot
rank constructed by Lin--Su \cite{LinSu} admit indecomposable lattice
lifts whose generic fibres are decomposable.

\subsection{Counterexamples to two GLS conjectures}
\label{subsec:gls-counterexamples}
We exhibit an example with connected Cartan matrix and minimal
symmetrizer that contradicts two GLS conjectures. The first,
\cite[p.~23, Conjecture~5.1]{GLSII}, predicts that
$\Rep(H,D\mathbf r)$ is irreducible for every rank vector
$\mathbf r$, equivalently that the locally free locus is dense.
The conjectural assertion in \cite[p.~23, Section~5.4]{GLSII}
predicts that, for an $H$-Schur root $\mathbf r$, the endomorphism
algebra of a general locally free $H$-module of rank $\mathbf r$ is
isomorphic to some vertex algebra $H_i$. Here $H$-Schur means that
indecomposable modules are dense in
$\Rep_{\mathrm{l.f.}}(H,\mathbf r)$, as in
\cite[Section~3.3]{GLSII}.

\begin{proposition}[{Counterexamples to two GLS conjectures}]
\label{prop:gls-counterexamples}\label{ex:gls-counterexamples}
Let
\[
 C^\dagger=\begin{pmatrix}2&-2&0\\-3&2&-2\\0&-2&2\end{pmatrix},
 \qquad D^\dagger=\op{diag}(3,2,2),
 \qquad 1\longrightarrow2\rightrightarrows3,
\]
and let $H=H(C^\dagger,D^\dagger,\Omega)$ for the displayed
orientation. Then the following hold.
\begin{enumerate}
\item The representation variety $\Rep(H,(3,2,2))$ has at least two
distinct irreducible components of dimension $20$. One is the
closure of the locally free rank-$(1,1,1)$ locus; another has a
dense open subset on which the loops at vertices $2$ and $3$
vanish.
\item The general locally free $H$-module of rank $(1,1,1)$ has
endomorphism algebra $k$. In particular, $(1,1,1)$ is an
$H$-Schur root, while its generic endomorphism algebra is not
isomorphic to any vertex algebra $H_i$.
\end{enumerate}
\end{proposition}
\begin{proof}
The equalities $2d_1=3d_2$ and $d_2=d_3$ show that
$D^\dagger$ is the minimal symmetrizer, and the datum is connected.
In column coordinates write the arrows as
$A_{12}:k^3\to k^2$ and $B_0,B_1:k^2\to k^2$.
The loop matrices $X_i$ satisfy
$X_1^3=0,\ X_2^2=X_3^2=0$ and $X_3B_j=B_jX_2\ (j=0,1)$.
The first-edge relation $X_2^2A_{12}=A_{12}X_1^3$ is redundant,
so $A_{12}$ is arbitrary.

Let $U_{\mathrm{nf}}$ be the locus where $B_0$ is invertible
and $B_0^{-1}B_1$ has two distinct eigenvalues. It is nonempty
and open in $\Rep(H,(3,2,2))$. On this locus, $X_2$ commutes
with $B_0^{-1}B_1$. A nilpotent matrix commuting with a matrix
with distinct eigenvalues is zero, so $X_2=0$ and then
$X_3=B_0X_2B_0^{-1}=0$. Conversely these equalities impose no
further restrictions on $X_1$, $A_{12}$, or the two arrows.
Thus $U_{\mathrm{nf}}$ is the product of the nilpotent cone in
$\mathfrak{gl}_3$, the space of $2\times3$ matrices, and a
nonempty open subset of pairs of $2\times2$ matrices. It is
irreducible, of dimension $6+6+8=20$.

Let $U_{\mathrm{lf}}=\Rep_{\mathrm{l.f.}}(H,(1,1,1))$.
This is the nonempty open locus on which all three loops are
regular nilpotent. The product of their conjugacy orbits is
irreducible of dimension $6+2+2$. Over any point of this product,
$A_{12}$ varies freely in dimension $6$, and each intertwining
space for $B_j$ has dimension $2$. These spaces form a vector
bundle over the product of the loop orbits. Hence
$U_{\mathrm{lf}}$ is irreducible and
$\dim U_{\mathrm{lf}}=(6+2+2)+(6+2+2)=20$.
The nonempty irreducible open subsets $U_{\mathrm{nf}}$ and
$U_{\mathrm{lf}}$ are disjoint. Their closures are therefore
distinct irreducible components, both of dimension $20$.
This proves (1), and refutes
\cite[p.~23, Conjecture~5.1]{GLSII} at the
$D^\dagger$-divisible dimension vector $(3,2,2)$.

For (2), let $J_m$ be the lower nilpotent Jordan block of size
$m$, and consider the locally free module $M$ given by
\[
 X_1=J_3,\qquad X_2=X_3=J_2,\qquad
 A_{12}=\begin{pmatrix}0&1&0\\0&0&1\end{pmatrix},\qquad
 B_0=I_2,\quad B_1=0.
\]
An endomorphism commuting with the loops has vertex matrices
\[
 P_1=x_0I_3+x_1J_3+x_2J_3^2,
 \qquad P_2=y_0I_2+y_1J_2.
\]
The arrow $B_0$ forces $P_3=P_2$, and the remaining condition is
\[
 P_2A_{12}-A_{12}P_1
 =\begin{pmatrix}
   -x_1&y_0-x_0&0\\
   -x_2&y_1-x_1&y_0-x_0
  \end{pmatrix}=0.
\]
It follows that $x_1=x_2=y_1=0$ and $y_0=x_0$, so
$\End_H(M)=k$. Since $U_{\mathrm{lf}}$ is irreducible,
upper semicontinuity of the endomorphism dimension gives a
dense open subset of $U_{\mathrm{lf}}$ on which every
endomorphism algebra is $k$. Such modules are indecomposable,
so $(1,1,1)$ is an $H$-Schur root. The vertex algebras are
$k[x]/(x^3)$ and $k[x]/(x^2)$, neither of which is isomorphic to
$k$. This refutes the prediction in
\cite[p.~23, Section~5.4]{GLSII}.
\end{proof}

The Schur-endomorphism counterexample already occurs on the first
two vertices. The third vertex and its two parallel arrows supply
the additional component needed to disprove irreducibility.

\begin{remark}[The generic-species comparison]
For the same datum, the corresponding standard generic species
has vertex fields $K_3,K_2,K_2$. Consider a rank-$(1,1,1)$ representation
whose first arrow and at least one parallel arrow are nonzero.
Every endomorphism is a tuple of vertex scalars. If one scalar
vanishes, compatibility with the nonzero arrow maps forces all
three to vanish. Thus every projection of the endomorphism
algebra to a vertex field is injective. This algebra is therefore
a finite field extension of $K$ embedded in each of $K_3,K_2,K_2$.
Its degree divides both $3$ and $2$, and consequently it is $K$.
Hence even on the generic fibre, a Schur endomorphism field need
not be one of the vertex fields.
\end{remark}

\section*{Acknowledgment}
The author thanks ChatGPT for suggesting the use of BPS cohomology
in the root multiplicity problem and for providing most of the
arguments in Section~\ref{sec:allrank-cohomology}.
The author also thanks Zengqiang Lin for sharing the example in
\cite{LinSu} before it appeared on arXiv.

\end{document}